\documentclass[a4paper]{amsart}
\input{commandoos.sty}

\begin{document}
	\author{A. Eggink}
	\title{Diophantine Maps} 
	\begin{abstract}
		To prove that Hilbert's tenth problem over a ring $R$ has a negative answer, usually the integers or another ring for which Hilbert's tenth problem has a negative solution is modelled inside the ring of interest. In this paper, we formalize this practice by introducing the notions of a Diophantine map and a Diophantine equivalence map. We compare the Diophantine case to the recursive case. We formalise a general version of Hilbert's tenth problem and show that we can transfer a positive or negative answer to Hilbert's tenth problem using effective Diophantine maps.  
	\end{abstract}
	
	\maketitle 
		
		\section{Introduction}
		
		\subsection{Hilbert's Tenth Problem}\phantom{=}\\
		Hilbert's tenth problem is the tenth problem on a list of problems that David Hilbert posed in 1900. In modern language, it asks for an algorithm to decide whether or not a Diophantine equation with coefficients in the integers has a solution in the integers. In 1970 it was proven by Yuri Matiyasevich, using the work of Julia Robinson, Martin Davis and Hilary Putnam, that such an algorithm can not exist and Hilbert's tenth problem has a negative answer. A full explanation of this proof can be found in \cite{H10PMurtyFodden}. In the proof different notions are related. The first one is Diophantine sets.
		\begin{defi} \th\label{introDefDiophSet} Let $(R,\L)$ be a set and a language (with an interpretation in $R$). A \textbf{term} of $\L$ consists of constants of $\L$ and variables, combined with functions of $\L$. A \textbf{basic formula} is of the form $(t_1,\ldots,t_k)\in R'$ where $R' \in \L$ is an $k$-ary relation (which can be equality) and the $t_i$ are terms of $\L$.\\
			A set $S \subset R^k$ is \textbf{Diophantine} if there exist basic formulas $f_1,\ldots,f_r$ such that
			$$S = \{(x_1,\ldots,x_k) \in R^k\mid \exists y_1,\ldots,y_l \in R\ (f_1\wedge f_2\wedge \ldots \wedge f_r)\}.$$   
		\end{defi}
		To recover the more standard definition of a Diophantine set, we take $R$ to be a ring and $\L$ to be $\L_R = \{0,1,+,\cdot\}$, the language of rings.\\
		The other notions are those of recursively enumerable sets and recursive sets. Recursively enumerable sets are sets that can be enumerated by an algorithm. Recursive sets are decidable sets, so sets for which there exists an algorithm that decides whether an element belongs to that set or not. We also have the notion of a recursive (or computable) function, which is informally an algorithm. For the precise definition and an introduction see appendix A of \cite{Shlapentokh}. Note that each Diophantine set is recursively enumerable, as we can just try all options $(\vec{x},\vec{y})$ one by one and output $\vec{x}$ if the equations are satisfied.\\
		The MRDP theorem, named after the four main contributors states:
		\begin{thm}[MRDP]
			Let $S \subset \Z^k$ be a recursively enumerable set. If we give $\Z$ the language of rings, then $S$ is a Diophantine set.
		\end{thm}
		Since there exists a recursively enumerable set that is not recursive, this implies that Hilbert's tenth problem over $\Z$ has a negative answer.
		
		\subsection{Diophantine Maps} \label{paragraphIntroDiophMap} \phantom{=}\\
		The next question then is: What if we replace $\Z$ by another ring? We need to be a bit careful with this since algorithms can not handle an uncountable amount of options. If our ring $R$ and corresponding language $\L$ behave well (more precisely, as in \th\ref{H10PexistenceH10P}), we say that Hilbert's tenth problem exists over $(R,\L)$.\\
		Hilbert's tenth problem over other rings has been addressed before. The main open problem in this area is Hilbert's tenth problem over $\Q$. It is known that Hilbert's tenth problem has a negative answer over $\F_q(T)$ by work of Pheidas \cite{PheidasFpT} and Videla \cite{H10PF2(T)}, over $\C(t_1,t_2)$ by Kim and Roush \cite{KRC(t1t2)} and over $R[T]$ for a commutative domain $R$ by Denef \cite{DenefRTchar0}, \cite{DenefRTcharp} and there are many more results. The common thing between these articles is that they involve constructing $\Z$, or another ring for which Hilbert's tenth problem has a negative solution, inside the ring of interest. \\
		Setting up a framework for this leads to my definition of a Diophantine map:
		\begin{defi} \th\label{introDefDiophMap} 
			Let $(R_1,\L_1)$ and $(R_2,\L_2)$ be sets with a language. A map $d: R_1 \ra R_2$ is called a \textbf{Diophantine map} if for all Diophantine sets $S \subset R_1^k$ the set $d(S) \subset R_2^k$ is Diophantine.
		\end{defi}
		Since this definition is hard to check, we have one of the main theorems:
		\begin{thm} \th\label{introDiophMap}
			Let $(R_1,\L_1)$ and $(R_2,\L_2)$ be sets with a language. Let $d: R_1 \ra R_2$ be an injective map such that:
			\begin{itemize}
				\item For all constants $c \in \L_1$, $d(c)$ can be defined by $\L_2$;
				\item For all functions $f: R_1^k \ra R_1$ of $\L_1$ the set $\{(d(x_1),\ldots,d(x_k),d(f(x_1,\ldots,x_k))) \mid \vec{x} \in R_1^k\}$ is Diophantine;
				\item For all relations $S \subset R_1^k$ of $\L_1$ the set $d(S)$ is Diophantine;
				\item The set $d(R_1) \subset R_2$ is Diophantine.
			\end{itemize}
			Then $d$ is a Diophantine map.
		\end{thm}
		The following result tells us that we can transfer a negative answer to Hilbert's tenth problem. 
		\begin{lemma} \th\label{introH10Pnegtransf} 
			Let $(R_1,\L_1)$ and $(R_2,\L_2)$ be sets with a recursive interpretation of a language. Let $d:R_1 \ra R_2$ be an injective Diophantine and recursive map. If $S \subset R_1$ is a Diophantine set that is not a recursive set, then $d(S)$ is also a Diophantine set that is not a recursive set. If furthermore for all $x \in d(S)$ the set $\{x\}$ is Diophantine, then Hilbert's tenth problem over $(R_2,\L_2)$ has a negative answer. 
		\end{lemma}  
		Since sometimes it is necessary to model elements in the first ring by sets in the second ring, this leads to the definition of a Diophantine equivalence map.
		\begin{defi} \th\label{introDefDiophEquiv} 
			Let $(R_1,\L_1)$ and $(R_2,\L_2)$ be sets with a language. Let $\sim$ be an equivalence relation on $R_2$ and let $\pi: R_2 \ra R_2\modsim$ be the projection. If $d: R_1 \ra R_2\modsim$ is a map such that for all $S \subset R_1^k$ Diophantine the set 
			$$\ol{d}(S) = \{\vec{x} \in R_2^k \mid \va 1\leq i \leq k,\ \pi(x_i) \in d(S)\}$$
			is a Diophantine set, then $d$ is a \textbf{Diophantine equivalence map}. 
		\end{defi}
		Just as before there is a set of conditions for a map to be a Diophantine equivalence map.
		\begin{thm} \th\label{introDiophEquiv} 
			Let $(R_1,\L_1)$ and $(R_2,\L_2)$ be sets with a language. Let $\sim$ be an equivalence relation on $R_2$ and let $\pi:R_2 \ra R_2\modsim$ be the projection. Let $d: R_1 \ra R_2\modsim$ be an injective map such that
			\begin{itemize}
				\item for all constants $c \in \L$, $\ol{d}(c) \subset R_2$ is a Diophantine set;
				\item for all $k$-ary functions $f \in \L$, $\ol{d}(\{(\vec{x},f(\vec{x}))\mid \vec{x} \in R_1^k\})$ is a Diophantine set;
				\item for all relations $S \subset \L$, $\ol{d}(S)$ is a Diophantine set;
				\item the set $\ol{d}(R_1)$ is Diophantine;
				\item the set $\{(x,y) \in R_2^2\mid x \sim  y\}$ or the set $\{(x,y)\in R_2^2 \mid x,y \in \ol{d}(R_1),\ x\sim y\}$ is Diophantine.
			\end{itemize}
			Then $d$ is a Diophantine equivalence map.
		\end{thm} 
		As expected we can also transfer a negative answer to Hilbert's tenth problem with a Diophantine equivalence map. For the exact statement see \th\ref{H10Pddiophequiv+rec+R1und=>R2und}.\\ 
		In order to transfer a positive answer to Hilbert's tenth problem, we need a way to construct the Diophantine equations of $d(S)$ from those of $S$. If this can be done, we call a Diophantine map effective. For the exact definition, see \th\ref{H10Pdefeffdioph}. We have the following theorem, which tells us that our main tool already gives us effective Diophantine maps.
		\begin{thm} \th\label{introEffDiophMap} 
			Take $(R_1,\L_1)$ and $(R_2,\L_2)$ sets with a finite language. Suppose that Hilbert's tenth problem exists over $(R_1,\L_1)$ and $(R_2,\L_2)$. Let $d: (R_1,\L_1) \ra (R_2,\L_2)$ be an injective map. 
			Suppose that for all the following sets:
			\begin{itemize}
				\item For all constants $c \in \L_1$, the set $\{d(c)\}$;
				\item For all functions $f: R_1^k \ra R_1$ of $\L_1$ the set $\{(d(x_1),\ldots,d(x_k),d(f(x_1,\ldots,x_k))) \mid \vec{x} \in R_1^k\}$;
				\item For all relations $S \subset R_1^k$ of $\L_1$ the set $d(S)$;
				\item The set $d(R_1)$,
			\end{itemize}
			an explicit description of the form 
			$$\{\vec{x} \in R_2^k \mid \ex \vec{y} \in R_2^l,\ \va 1\leq i \leq n,\ (t_{i1}(\vec{x},\vec{y}),\ldots,t_{in_i}(\vec{x},\vec{y})) \in S_i\}$$
			with $t_{ij}$ terms of $\L_2$ and $S_i$ relations of $\L_2$ is known.\\
			Then the map $d$ is an effective Diophantine map.
		\end{thm} 
		The following result shows that we can transfer a positive answer to Hilbert's tenth problem backwards using these effective Diophantine maps:
		\begin{thm} \th\label{introH10Ppostransf} 
			Let $(R_1,\L_1)$ and $(R_2,\L_2)$ be sets with a countable language. Suppose that Hilbert's tenth problem exists over $(R_1,\L_1)$ and $(R_2,\L_2)$. If Hilbert's tenth problem over $(R_2,\L_2)$ has a positive answer and $d:R_1 \ra R_2$ is an effective Diophantine map, then Hilbert's tenth problem over $(R_1,\L_1)$ has a positive answer.
		\end{thm}
		
		\begin{remark}
			This theorem tells us that effective Diophantine maps play the role of ``morphisms" if we look at the decidability of positive existential theories.
		\end{remark} 

		Besides the question of Hilbert's tenth problem over a certain ring $R$, we can also ask whether the recursively enumerable sets of $R$ are the same as the Diophantine sets. We say that MRDP holds in $R$ if this is the case.\\
		Since the MRDP theorem connects Diophantine sets with recursively enumerable sets, it can be asked if we can connect Diophantine maps to recursive functions. To answer that, we need two more notions.
		\begin{defi} \th\label{introDefREpres} 
			We call a function $f: R_1 \ra R_2$ \textbf{recursively enumerable preserving}, or \textbf{r.e. preserving} for short, if for all recursively enumerable sets $S \subset R_1$, the set $f(S) \subset R_2$ is recursively enumerable. 
		\end{defi}
		\begin{defi} \th\label{introDefGivenbyDiophExpr} Let $(R_1,\L_1)$ and $(R_2,\L_2)$ be sets with a language. Take maps $d_1: R_1 \ra R_2$ and $d_2: R_1^k \ra R_2^l$. We say that $d_2$ is \textbf{given by a Diophantine expression} in $d_1$ if there exist basic $\L_2$-formulas $f_1,\ldots,f_n$ such that 
			$$d_2(\vec{x}) = \vec{z} \ \Leftrightarrow\ \ex \vec{y} \in R^{l'},\ \va 1\leq i \leq n,\ f_i(d_1(\vec{x}),\vec{y},\vec{z})$$ 
			If $k=l=1$, we will denote this relation with $d_2 \geq_D d_1$ or $d_1 \leq_D d_2$.
		\end{defi}
		The relations between these notions can then be shown in a diagram. Let $(R_1,\L_1)$ and $(R_2,\L_2)$ be infinite countable rings and let $\phi:R_1 \ra R_2$ be a bijective recursive function. (If $R_1=R_2$ we can take the identity, but such a $\phi$ always exists by \th\ref{recrecbij}.) Then we get:\\
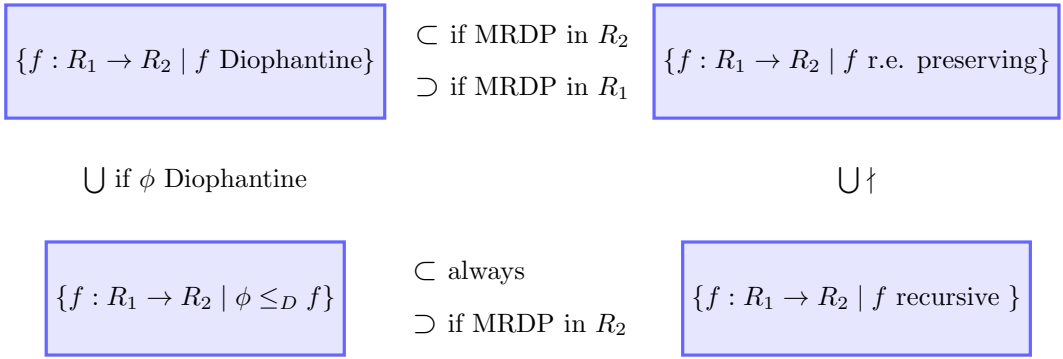
		\begin{center} 
			\begin{tikzpicture}[squarednode/.style={rectangle, draw=blue!60, fill=blue!10, very thick, minimum size=15mm}]
				
				\node[squarednode] (Dioph) {$\{ f: R_1 \ra R_2\mid f$ Diophantine$\}$} ;
				\node[squarednode] (REpres) [right=3.5cm of Dioph] {$\{ f: R_1 \ra R_2\mid f$ r.e. preserving$\}$};
				
				\node (DiophREpres) [above right=-0.7 cm and 0.25cm of Dioph] {\Large $\subset$ \normalsize if MRDP in $R_2$};
				\node (DiophREpres2) [below right=-0.7 cm and 0.25cm of Dioph] {\Large $\supset$ \normalsize if MRDP in $R_1$};
				
				\node (RecREpres) [below=0.5cm of REpres] {$\bigcup\nmid$};
				\node (DiophDiophexpr) [below=0.5 of Dioph] {$\bigcup$ if $\phi$ Diophantine};
				
				\node[squarednode] (Diophexpr) [below=0.5cm of DiophDiophexpr] {$\{ f: R_1 \ra R_2\mid \phi \leq_D f \}$};
				\node[squarednode] (recursive) [below=0.5cm of RecREpres] {$\{ f: R_1 \ra R_2\mid f$ recursive $\}$};
				
				\node (DiophexprRec) [above right=-0.7 cm and 0.75cm of Diophexpr] {\Large$\subset$ \normalsize always};
				\node (DiophexprRec2) [below right=-0.7 cm and 0.75cm of Diophexpr] {\Large$\supset$ \normalsize if MRDP in $R_2$};
			\end{tikzpicture}
		\end{center}
		\normalsize 
		
		\section{Acknowledgements} 
		I thank my supervisor, Gunther Cornelissen, for suggesting the topic, his supervision during the project and his help in writing this paper. \\
		I also thank Philip Dittmann for his question about the need of compatibility between the ring structure of $R_0$ and the structure of $(R,\L)$ in the existence of Hilbert's tenth problem. This led to \th\ref{H10Pdefcompatibility} and \th\ref{H10PchoicerecreprR0}.
	
	\section{Conventions}
	I often use commas instead of wedges to denote and-statements, for example $\{(x,y) \in \Z^2 \mid x>0,\ y>0\}$ is the set of $(x,y)$ for which $x>0$ and $y>5$. \\
	I assume that a ring has a $1$, but it does not have to be commutative.\\
	By assumption each language contains a constant $0$ and a relation 
	``$=$", where often ``$=$" is left out from the notation.
	
	\section{Notation} 
	\begin{itemize}
		\item $\N = \Z_{\geq 0}$;
		\item $\L_R = \{0,1,+,\cdot\}$ is the language of rings;
		\item $c^R,f^R$ or $S^R$ is the interpretation of a constant $c$, function $f$ or relation $S$ of a language $\L$ in a set $R$;
		\item $S^c$ is the complement of $S$;
		\item $p_n$ is the $n$-th prime number;
		\item $(\mu x)(R(x))$ is the minimizing operator, which gives the least $x$ such that the relation $R(x)$ is true;
		\item $\rho(i,n)$ gives the exponent of $p_i$ in the prime factorization of $n$.
		\item $S_c(\L)$ is the set of constants of a language $\L$;
		\item $S_c'(\L) = \{x\in R\mid \{x\}\text{ is a Diophantine set}\}$, for a language $\L$ with an interpretation in the set $R$; 
		\item $\Frac(R)$ is the field of fractions of $R$;
		\item $\ol{d}(S) = \{\vec{x} \in R_2^k \mid (\pi(x_1),\ldots,\pi(x_k)) \in d(S)\}$, for a Diophantine equivalence map $d:R_1 \ra R_2\modsim$ and a set $S \subset R_1^k$;
		\item $t[x/y]$ for $t$ a term, $x$ a variable and $y$ a variable or an element of $S_c'(\L)$ is the term where each occurrence of $x$ is replaced with a $y$.
	\end{itemize}
	
	\section{Definitions and Properties of Diophantine Sets and Maps} 
	
	We start this section by introducing Diophantine sets and proving their basic properties. After that we define Diophantine maps and give some basic properties. Furthermore, we prove a theorem that is our main tool for getting Diophantine maps. Since the class of Diophantine maps is sometimes too broad, we introduce the notion of ``a map given by a Diophantine expression in another map". Then we move on to Diophantine equivalence maps, a generalization of Diophantine maps that sends an element to an equivalence class for some equivalence relation. Finally, we introduce Diophantine models and equivalences of Diophantine models. A lot of the proofs of lemmas in this section and some in further sections are straightforward, so we omit them. 
	
	\subsection{Diophantine Sets}\phantom{=}\\
	Recall the general definition of a Diophantine set from the introduction, \th\ref{introDefDiophSet}. 
	\begin{defi} \th\label{diophDefDiophSet} Let $(R,\L)$ be a set and a language (with an interpretation in $R$). A \textbf{term} of $\L$ consists of constants of $\L$ and variables, combined with functions of $\L$. A \textbf{basic formula} is of the form $(t_1,\ldots,t_k)\in R'$ where $R' \in \L$ is an $k$-ary relation (which can be equality) and the $t_i$ are terms of $\L$.\\
		A set $S \subset R^k$ is Diophantine if there exist basic formulas $f_1,\ldots,f_r$ such that
		$$S = \{(x_1,\ldots,x_k) \in R^k\mid \exists y_1,\ldots,y_l \in R\ (f_1\wedge f_2\wedge \ldots \wedge f_r)\}.$$   
	\end{defi}
	\begin{example} 
		Take $R$ a ring and $\L = S_c\cup \{+,\cdot,=\}$, where $S_c$ is a set of constants containing $0$ and $1$. Then the terms $t_i$ are polynomials with coefficients in $S_c$ and the $f_i$ are of the form $t_{i1} = t_{i2}$, so equivalently $t_{i1}-t_{i2}=0$. In this case we get the more standard definition of a Diophantine set:
		$$S = \{\vec{x} \in R^k\mid \ex \vec{y} \in R^{k'}\text{ s.t. }\va i\ P_i(\vec{x},\vec{y})=0\},$$
		where the $P_i$ are polynomials with coefficients in $S_c$ and terms may occur multiple times. 
	\end{example}
	
	\begin{assumption}
		We will assume that all languages contain a constant $0$ and the relation ``$=$". We will omit ``$=$" if we specify the symbols of a language.\\ 
		If on a ring no language is specified, it will have the language of rings $\L_R = \{0,1,+,\cdot\}$.
	\end{assumption}
	
	\begin{remark}
		Note that in the main applications, $R$ will be a ring and $\L$ will contain the language of rings, but this is not necessary. The only constraints are that $R$ is a set and that $\L$ contains a constant $0$ and the (from notation omitted) relation ``$=$''. 
	\end{remark}
	
	\begin{notation}
		Let $R$ be a set with an interpretation of a language $\L$. We denote with $S_c = S_c(\L)$ the set of constants of $\L$. We define $S_c' = S_c'(\L)$ by setting $x \in S_c'(\L)$ if $\{x\} \subset R$ is a Diophantine set.
	\end{notation}
	
	\begin{remark}
		All elements of $S_c'(\L)$ can be given by one or more Diophantine formulas, so we can use them as constants in our terms without changing which sets are Diophantine.
	\end{remark}
	
	\begin{example} 
		It can be useful to work with extra relations besides equality.
		A simple example is $\Z$ with language $\L = \{0,1,+,|\}$, where $|$ is the divisor relation, so $n|m$ if and only if $\ex\ k, nk=m$.
		We can also take a prime number $p \in \Z_{\geq 0}$ and define the relation $|_p$ by $n|_pm$ if and only if $\ex s \in \Z_{\geq 0}\ np^s=m$. Then we can look at $\N = \Z_{\geq 0}$ with language $\L'= \{0,1,+,|_p\}$. This structure is used to prove that Hilbert's tenth problem over $\F_{p^n}(t)$ has a negative answer, see \cite{PheidasN} and \cite{PheidasFpT}.  
	\end{example}
	
	\begin{prop} \th\label{diophsystem=single}
		Let $R$ be an integral domain and let $K$ be its field of fractions. Suppose that $K$ is not algebraically closed. Then for all $P_1,\ldots,P_n \in R[\vec{x}]$ there exist a polynomial $Q \in R[\vec{x}]$ such that $P_1(\vec{x}) = \ldots = P_n(\vec{x}) = 0$ if and only if $Q(\vec{x}) = 0$. 
	\end{prop}
	
	\begin{proof}
		Since $K$ is not algebraically closed, there exists a non-constant polynomial $f \in K[X]$ such that $f$ has no zeros in $K$. Suppose $f= \sum_{i=0}^k a_iX^i$, with $a_k\not=0$ and let $h(X,Y) = \sum_{i=0}^k a_iX^iY^{k-i}$ be the homogenous version of $f$. Suppose $h(a,b)=0$. If $b\not=0$, we have that $b$ is invertible. This gives $0=b^{-k}h(a,b) = \sum_{i=0}^k a_i\left(\frac{a}{b}\right)^i = f(ab^{-1})$. But $ab^{-1} \in K$, contradicting the fact that $f$ has no zeros in $K$. If $b=0$, we get $a_ka^k=0$, so $a=0$. This implies that $(0,0)$ is the only zero of $h$. \\
		Now we get that $P_1(\vec{x}) = P_2(\vec{x}) = 0$ is equivalent to $h(P_1(\vec{x}),P_2(\vec{x})) = 0$ and this is a single polynomial equation. By induction on $n$ there exists a polynomial $f \in K[\vec{x}]$ such that $P_1(\vec{x}) = \ldots = P_n(\vec{x}) = 0$ if and only if $f(\vec{x}) = 0$ for all $\vec{x} \in K^k \supset R^k$. Multiplying $f$ by a nonzero factor $b$ to clear the denominators gives $bf \in R[x]$ and as $R$ is an integral domain $bf(\vec{x}) = 0$ if and only if $f(\vec{x}) = 0$. So taking $bf = Q$ gives the result.
	\end{proof}
	
	\begin{cor} \th\label{diophsystem=singlecor}
		Let $R$ be an integral domain and $\L$ the language of rings with set of constants $S_c$. A system of polynomials with coefficients in $S_c'$ is equivalent to a single polynomial with coefficients in $S_c'$ if the algebraic closure of $\Frac(S_c')$ is not a subset of $\Frac(R)$.
	\end{cor}
	
	\begin{proof}
		We use the approach of the proof of \th\ref{diophsystem=single}, with $K = \Frac(S_c')$ and note that we should have no zeros inside $\Frac(R)$.
	\end{proof} 
	
	
	\begin{lemma} \th\label{diophintersection} Let $(R,\L)$ be a set with a language. If $S_1,S_2 \subset R^k$ are Diophantine sets, then $S_1 \cap S_2$ is Diophantine. \qed
	\end{lemma}
	
	\begin{defi}
		Let $(R,\L)$ be a set with a language. We say that the $\L$-formulas $f_1(\vec{x})$ and $f_2(\vec{x})$ are \textbf{equivalent} over $R$ if
		$$\{\vec{x} \in R^k\mid f_1(\vec{x})\text{ is true}\} = \{\vec{x} \in R^k\mid f_2(\vec{x})\text{ is true}\}.$$
	\end{defi}
	
	\begin{defi}
		Let $(R,\L)$ be a set with a language. We say that $\L$ is \textbf{one-sided} over $R$ if for all Diophantine formulas of $(R,\L)$ there exist an equivalent Diophantine formula of the form
		$$\ex \vec{y} \in R^k,\ t_1(\vec{x},\vec{y}) = 0,\ldots, t_n(\vec{x},\vec{y}) = 0,$$
		where the $t_i$ are $\L$-terms. 
	\end{defi}
	
	\begin{lemma} \th\label{diophaddinv=>one-sided} 
		Let $(R,\L)$ be a set with a language such that $\L$ contains no relations other then equality. Suppose that $+ \in \L$ and that $R$ contains additive inverses. Then $\L$ is one-sided over $R$. \qed 
	\end{lemma}
	
	\begin{lemma} \th\label{diophunion}
		Let $(R,\L)$ be a set with a language such that $\L$ is one-sided over $R$. Suppose that $\cdot \in \L$ and $a\cdot b = 0$ implies $a=0$ or $b=0$. (If we work over a ring, this is equivalent to $R$ being an integral domain.) Let $S_1,S_2 \subset R^k$ be Diophantine sets. Then $S_1 \cup S_2$ is Diophantine. \qed 
	\end{lemma}
	
	\begin{lemma} \th\label{diophfiniteset} Let $\L$ be the language of rings with set of constants $S_c$. Let $R$ be an integral domain and let $S \subset (S_c')^k$ be finite. Then $S$ is a Diophantine set. \qed 
	\end{lemma}
	
	\begin{lemma} \th\label{diophcartRR}
		Let $(R,\L)$ be a set with a language. If $S_1 \subset R^k$ and $S_2 \subset R^{k'}$ are Diophantine, then $S_1 \ti S_2 \subset R^{k+k'}$ is Diophantine. \qed 
	\end{lemma}
	
	\begin{lemma} \th\label{diophprojectionRR} Let $S \subset R^{k} \ti R^l$ be a Diophantine set. Then the projections $\pi_1(S) \subset R^k$ and $\pi_2(S) \subset R^l$ are Diophantine sets. \qed 
	\end{lemma} 
	
	\begin{remark} There is a big difference between sets $S \subset R^k$ that are Diophantine if we give $R^k$ the language of rings and sets that are Diophantine if we give the language of rings to $R$. In the first case we only have "diagonal" functions and constants, so all terms are of the form $(t,\ldots,t)$. In the second case, we have a lot more freedom, as we may use $x_i \in R$ in our terms, instead of only the vector $\vec{x} \in R^k$. \\
		For example, let $k=2$. In the second case $\{(0,1)\} = \{(x,y) \in R^2\mid x=y-1=0\}$ is a Diophantine set. In the first case using the term $x$ will give $\{\vec{x} \in R^2\mid \vec{x} = \vec{0}\}$, so we get the Diophantine set $\{(0,0)\}$. For $R^2$ with the language of rings $\{(0,1)\}$ is in fact not Diophantine.
	\end{remark}
	
	\begin{remark} The previous remark shows that we might need another language on $R_1 \ti R_2$ then just the language of rings, since then for $S_1 \subset R_1$ and $S_2 \subset R_2$ Diophantine sets, we do not in general have that $S_1 \ti S_2 \subset R_1 \ti R_2$ is Diophantine. 
		Suppose we define a language on $R_1 \ti R_2$. If we then choose $R_1$ equal to $R_2$, we would want that the Diophantine sets defined by this language are the same as the one defined by $(R_1,\L_R)$. This is not possible, since we cannot define functions that take elements of $R_1$ and $R_2$. However we can define a language $\L^{(k)}$ on $R^k$, using the language $\L$ on $R$, such that the Diophantine sets under $\L^{(k)}$ are the same as the Diophantine sets induced by $\L$.
	\end{remark} 
	
	\begin{notation}
		For a constant $c$, function $f$ or relation $S$ in a language $\L$, we denote the interpretation in a set $R$ as $c^R$, $f^R$ respectively $S^R$.
	\end{notation}
	
	\begin{defi} 
		Let $(R,\L)$ be a set with an interpretation of the language $\L$. We define the language $\L^{(k)}$ and interpret it on $R^k$ as follows:
		\begin{itemize}
			\item For all $c \in \L$ a constant, $c \in \L^{(k)}$ is a constant and $c^{R^k} = (c^R,0^R,\ldots,0^R)$;
			\item For all $f \in \L$ an $r$-ary function, $f \in \L^{(k)}$ is an $r$-ary function and \\$f^{R^k}(\vec{x_1},\ldots,\vec{x_k}) = (f^R(x_{11},\ldots,x_{k1}),0^R,\ldots,0^R)$;
			\item For all $S \in \L$ an $r$-ary relation, $S \in \L^{(k)}$ is an $r$-ary relation and \\$S^{R^k} = \{((x_1,\vec{0}),\ldots,(x_r,\vec{0})) \in R^{kr} \mid \vec{x} \in S^R\}$;
			\item For $1\leq i \leq k$ the $1$-ary function $\pi_i$ is in $\L^{(k)}$ and $\pi_i(\vec{x}) = (x_i,0^R,\ldots,0^R)$.
		\end{itemize}
	\end{defi}
	
	\begin{lemma} \th\label{diophRkRsamediophsets} 
		The Diophantine subsets of $(R^k)^l$ with the language $\L^{(k)}$ are the same as the Diophantine subsets of $R^{kl}$ with the language $\L$. 
	\end{lemma}
	
	\begin{proof}
		Let $S \subset R^{kl}$ be a Diophantine set for the language $\L$. Then there exist $\L$-terms $t_{ij}$ with $1 \leq i \leq n$, $1\leq j \leq m_i$ and relations $S_i$ of $\L$ (not necessary distinct) such that 
		\begin{align*} 
			S = \{(x_{11},\ldots,x_{1k},\ldots,x_{l1},\ldots,x_{lk})\in R^{kl}\mid \ex \vec{y} \in R^{l'},\ \va 1\leq i\leq n,\phantom{===}\\ (t_{i1}(x_{11},\ldots,x_{lk},\vec{y}),\ldots,t_{im_i}(x_{11},\ldots,x_{lk},\vec{y}))\in S_i^R\}.
		\end{align*} 
		If we replace $x_{j'i'}$ with $\pi_{i'}(\vec{x_{j'}})$ in all the terms $t_{ij}$ we get $\L$-terms $t_{ij}'$ such that for all $1 \leq i \leq n$, $1\leq j \leq m_i$ we have
		$$t_{ij}(x_{11},\ldots,x_{1k},\ldots,x_{l1},\ldots,x_{lk}) = t_{ij}'(\vec{x_1},\ldots,\vec{x_l}).$$
		Taking the vector $\vec{y} \in (R^k)^{l'}$, so only using the first coordinate, gives
		\begin{align*}
			S &= \{(\vec{x_1},\ldots,\vec{x_l})\in (R^k)^l \mid \ex \vec{y}\in (R^k)^{l'}, \va 1 \leq i \leq n,\\ &\phantom{====}((t_{i1}'(\vec{x_1},\ldots,\vec{x_l},\vec{y}),\vec{0}),\ldots,(t_{im_i}'(\vec{x_1},\ldots,\vec{x_l},\vec{y}),\vec{0})) \in S_i^{R^k} \},
		\end{align*} 
		which is a Diophantine set for the language $\L^{(k)}$.\\
		Let $S \subset R^{kl}$ be a Diophantine set for the language $\L^{(k)}$. Then there exist $\L^{(k)}$-terms $t_{ij}$ and $\L^{(k)}$-relations $S_i$  such that 
		$$S = \{(\vec{x_1},\ldots,\vec{x_l})\in (R^k)^l\mid\ \ex\vec{y} \in (R^k)^{l'},\ \va 1\leq i \leq n,\ (t_{i1}(\vec{x_1},\ldots,\vec{x_l},\vec{y}),\ldots,t_{im_i}(\vec{x_1},\ldots,\vec{x_l},\vec{y}))\in S_i\}.$$
		By the construction of $\L^{(k)}$, an $\L^{(k)}$-term $t_{ij}$ either contains no functions or constants, or the last operation in the construction of $t_{ij}$ is a function or constant. In the first case we get that $t_{ij}= \vec{x}$ and we can write $t_{ij} = (t_{ij}',x_2,\ldots,x_k)$ with $t_{ij}'=x_1$ an $\L$-term. In the second case we get $t_{ij} = (t_{ij}',0,\ldots,0)$ where $t_{ij}'$ is an $\L$-term, as we have used constants and function of $\L$ and coordinate wise variables, defined via $\pi_i$. By construction of $S_i^{R^k}$ we have that $(t_{i1},\ldots,t_{im_i})\in S_i^{R^k}$ if and only if $(t_{i1}',\ldots,t_{im_i}')\in S_i^{R}$ and possible some extra constraints $x_{ij}=0$. This is a set of basic $\L$-formulas, so $S$ is also a Diophantine set with language $\L$, as we can still take $\vec{y}\in (R^k)^l = R^{kl}$.   
	\end{proof}
	
	\begin{defi} Let $(R_1,\L_1)$ and $(R_2,\L_2)$ be sets with an interpretation of a language. Then $\L_1 \sqcup \L_2$, the disjoint union of $\L_1$ and $\L_2$, has an interpretation on $R_1\ti R_2$ by
		\begin{itemize}
			\item a constant $c \in \L_1$ is $(c^{R_1},0^{R_2})$;
			\item a constant $c \in \L_2$ is $(0^{R_1},c^{R_2})$;
			\item an $r$-ary function symbol $f \in \L_1$ is the $r$-ary function $(x_i,y_i)_{i=1}^r \ra (f^{R_1}(\vec{x}),0^{R_2})$;
			\item an $r$-ary function symbol $f \in \L_2$ is the $r$-ary function $(x_i,y_i)_{i=1}^r \ra (0^{R_1},f^{R_2}(\vec{y}))$;
			\item an $r$-ary relation symbol $S \in \L_1$ is the $r$-ary relation $S^{R_1} \ti \{0^{R_2}\}$;
			\item an $r$-ary relation symbol $S \in \L_2$ is the $r$-ary relation $\{0^{R_1}\}\ti S^{R_2}$.
		\end{itemize}
		What this interpretation of $\L_1\sqcup \L_2$ does is make a copy of $(R_1,\L_1)$ in the first coordinate, where the second coordinate is always zero and make a copy of $(R_2,\L_2)$ in the second coordinate, while the first coordinate is always zero.
	\end{defi}
	
	There seem to exist other natural definitions for a language on $R_1\ti R_2$. We will take a look at two of them in the next two remarks.
	\begin{remark} 
		We begin with the language $\L_1\ti \L_2$. Suppose that $\L_1$ and $\L_2$ do not contain any relations other than equality. A way to define and interpret this language could be
		\begin{itemize}
			\item For $c_i \in \L_i$ constants, we let $(c_1,c_2) = (c_1^{R_1},c_2^{R_2})$ be a constant;
			\item For $c_1 \in \L_1$ a constant and $f_2 \in \L_2$ an $r$-ary function we let $(c_1,f_2)$ be the $r$-ary function defined by $(c_1,f_2)((x_1,y_1),\ldots,(x_r,y_r)) = (c_1^{R_1},f_2^{R_2}(y_1,\ldots,y_r))$;
			\item For $f_1 \in \L_1$ an $r$-ary function and $c_2 \in \L_2$ a constant we let $(f_1,c_2)$ be the $r$-ary function defined by $(f_1,c_2)((x_1,y_1),\ldots,(x_r,y_r)) = (f_1^{R_1}(x_1,\ldots,x_r),c_2^{R_2})$;
			\item For $f_1 \in \L_1$ an $r_1$-ary function and $f_2 \in \L_2$ an $r_2$-ary function we let $(f_1,f_2)$ be an $r$-ary function, where $r = \max(r_1,r_2)$. We then define\\ $(f_1,f_2)((x_1,y_1),\ldots,(x_r,y_r)) = (f_1^{R_1}(x_1,\ldots,x_{r_1}),f_2^{R_2}(y_1,\ldots,y_{r_2}))$.
		\end{itemize}
		This definition is less natural than it seems. One would expect that the terms of this function are $(t_1,t_2)$ with $t_i$ an $\L_i$-term, but this is not true in general. If we want to create the term $(0,y)$, we want to say it is $(0,"id")(x,y)$ as only the pair $(x,y)$ is a variable and not $x$ and $y$ themselves, but we do not have an identity function. \\
		If $\L_1 = \L_2 = \L_R$ and $R_1$ and $R_2$ are rings, we can write $(0,y) = (0,"id")(x,y) = (0,+)((x,y),(0,0))$, so we create an identity function by using that $x+0 = x$. In this case we also have that the terms of $\L_1\ti \L_2$ are of the form $(t_1,t_2)$ with $t_i$ an $\L_i$-term, which can be proven by induction.
	\end{remark}
	
	\begin{remark}  
		The lack of the identity function in the previous language leads to another possible language on $R_1\ti R_2$, namely $(\L_1\sqcup \L_2)^*$. It contains no constants and the following functions and relations, which we interpret in $R_1 \ti R_2$:
		\begin{itemize}
			\item For each constant $c \in \L_1$ the $1$-ary function $c_1$ and $c_1^{R_1\ti R_2}:(x,y) \ra (c^{R_1},y)$;
			\item For each constant $c \in \L_2$ the $1$-ary function $c_2$ and $c_2^{R_1\ti R_2}:(x,y) \ra (x,c^{R_2})$;
			\item For each $r$-ary function $f \in \L_1$ the $r$-ary function $f_1$ and\\ $f_1^{R_1 \ti R_2}((x_1,y_1),\ldots,(x_r,y_r)) = (f^{R_1}(x_1,\ldots,x_r),y_1)$;
			\item For each $r$-ary function $f \in \L_2$ the $r$-ary function $f_2$ and\\ $f_2^{R_1 \ti R_2}((x_1,y_1),\ldots,(x_r,y_r)) = (x_1,f^{R_2}(y_1,\ldots,y_r))$;
			\item For each $r$-ary relation $S \in \L_1$ the $r$-ary relation $S$ and\\ $S^{R_1\ti R_2} = \{(x_1,y_1,\ldots,x_r,y_r) \in (R_1\ti R_2)^r \mid \vec{x} \in S^{R_1}, \va 1\leq i \leq r,\ y_i = 0^{R_2}\}$;
			\item For each $r$-ary relation $S \in \L_2$ the $r$-ary relation $S$ and\\ $S^{R_1\ti R_2} = \{(x_1,y_1,\ldots,x_r,y_r) \in (R_1\ti R_2)^r \mid \vec{y} \in S^{R_2}, \va 1\leq i \leq r,\ x_i = 0^{R_1}\}$.
		\end{itemize}
		This language has terms $(t_1,t_2)$ with $t_i$ an $\L_i$-term, as we can create $(t_1,y_1)$ by applying constants and functions of $\L_1$ to $(x_i,y_i)$ and then we can create $(t_1,t_2)$ by applying constants and functions of $\L_2$ to $(t_1,y_1),(x_2,y_2),\ldots,(x_n,y_n)$.\\
		If $S_i$ is a relation of $\L_1$, we have that $((t_1,s_1),\ldots,(t_n,s_n)) \in S_i$ is equivalent to $((t_1,0),\ldots,(t_n,0)) \in S_i$ and $(0,s_1)=\ldots = (0,s_n)= 0$. This is a set of Diophantine formulas in $\L_1\sqcup \L_2$. By symmetry, something similar holds for $S_i \in \L_2$, so it turns out that we only get more Diophantine expressions and not more Diophantine sets. As the terms of $\L_1 \sqcup \L_2$ are simpler than those of $(\L_1 \sqcup \L_2)^*$, we will use the first as our standard language on $R_1 \ti R_2$. 
	\end{remark}
	
	\begin{lemma} \th\label{diophcartsetR1R2} Let $(R_1,\L_1)$ and $(R_2,\L_2)$ be sets with languages. Equip $R_1 \ti R_2$ with $\L_1 \sqcup \L_2$. If $S_1 \subset R_1^k$ and $S_2 \subset R_2^k$ are Diophantine, then $S_1 \ti S_2 \subset (R_1\ti R_2)^k$ is a Diophantine set. \qed 
	\end{lemma}
	
	\begin{lemma} \th\label{diophprojectionR1R2} Let $(R_1,\L_1)$ and $(R_2,\L_2)$ be sets with a language such that $\emptyset \subset R_1$ and $\emptyset \subset R_2$ are both a Diophantine set or both not a Diophantine set. Equip $R_1\ti R_2$ with $\L_1 \sqcup \L_2$ and let $S \subset (R_1\ti R_2)^k$ be a Diophantine set. Then the projections $\pi_1(S) \subset R_1^k$ and $\pi_2(S) \subset R_2^k$ are Diophantine sets.
	\end{lemma} 
	
	\begin{proof} We have that the terms of $\L_1 \sqcup \L_2$ are of the form $(t,0)$ or $(0,s)$ where $t$ is an $\L_1$-term and $s$ is an $\L_2$-term. Furthermore the relations are of the form $T_i' = \{(x_1,y_1,\ldots,x_{k_i},y_{k_i})\mid \vec{x} \in T_i\}$ with $T_i$ a relation of $\L_1$ or $S_i' = \{(x_1,y_1,\ldots,x_{k_i},y_{k_i})\mid \vec{y} \in S_i\}$ with $S_i$ a relation of $\L_2$. This means we can write 
		\begin{align*} 
			S &= \{(x_i,y_i)_{i=1}^k \in (R_1\ti R_2)^k \mid \ex (u_i,v_i)_{i=1}^{k'} \in (R_1\ti R_2)^{k'},\\ 
			&\phantom{=====}\va 1\leq i \leq n,\ ((t_{i1},s_{i1}),\ldots,(t_{im_i},s_{im_i})) \in T_i'\ti \{0\},\\
			&\phantom{======}\va 1\leq i \leq n',\ ((t_{i1}',s_{i1}'),\ldots,(t_{im_i'}',s_{im_i'}')) \in \{0\}\ti S_i'\}
		\end{align*} 
		Suppose first that $S$ is not the empty set. The statement $((t_{i1},s_{i1}),\ldots,(t_{im_i},s_{im_i})) \in T_i'$ 
		is true if and only if $(t_{i1},\ldots,t_{im_i}) \in T_i$ and $s_{i1}=\ldots = s_{im_i} = 0$. There exist $\vec{y} \in R_2^k, \vec{w} \in R_2^{k'}$ such that $s_{i1}=\ldots=s_{im_i}=0$ is true for some $\vec{x} \in R_1^k$, so also for all $\vec{x} \in R_1^k$, since it has no occurrence of $\vec{x}$. In the same way we have that $((t_{i1}',s_{i1}'),\ldots,(t_{im_i'}',s_{im_i'}')) \in \{0\}\ti S_i'$ is true if and only if $(s_{i1}',\ldots,s_{im_i'}')\in S_i$ and $t_{i1}'=\ldots = t_{im_i'}'=0$. The statement $(s_{i1}',\ldots,s_{im_i'}')\in S_i$ contains no $\vec{x}$, so since $S$ is non-empty, it is true for all $\vec{x} \in R_1^k$. Combining this gives
		$$\pi_1(S) = \{\vec{x} \in R_1^k \mid \ex \vec{y} \in R_1^{k'},\ \va 1\leq i\leq n,\ (t_{i1},\ldots,t_{im_i}) \in T_i,\ \va 1\leq i \leq n',\ t_{i1}'=\ldots = t_{im_i'}'=0\},$$
		which is a Diophantine set in $(R_1,\L_1)$. \\
		If $S = \emptyset$, then $\emptyset \subset R_1$ or $\emptyset \subset R_2$ is Diophantine. This gives by assumption that $\emptyset \subset R_1$ and $\emptyset \subset R_2$ are Diophantine sets, so $\pi_1(S)=\pi_1(\emptyset) = \emptyset$ is a Diophantine set. \\
		It follows by symmetry that $\pi_2(S)$ is a Diophantine set.
	\end{proof}
	
	\begin{lemma} \th\label{diophcartsplit}
		Let $(R_1,\L_1)$ and $(R_2,\L_2)$ be sets with a language such that $\emptyset \subset R_1$ and $\emptyset \subset R_2$ are both a Diophantine set or both not a Diophantine set. Equip $R_1 \ti R_2$ with $\L_1 \sqcup \L_2$. If $S \subset (R_1\ti R_2)^k$ is a Diophantine set, then $S = S_1 \ti S_2$ with $S_i \subset R_i^k$ a Diophantine set.
	\end{lemma}
	
	\begin{proof}
		We have seen in the proof of the previous lemma that we can write 
		\begin{align*} 
			S &=  \{(x_i,y_i)_{i=1}^k \in (R_1\ti R_2)^k \mid \ex (u_i,v_i)_{i=1}^{k'} \in (R_1\ti R_2)^{k'},\\ 
			&\phantom{=====}\va 1\leq i \leq n,\ (t_{i1},\ldots,t_{im_i}) \in T_i,\ s_{i1}=\ldots = s_{im_i} = 0 \\
			&\phantom{======}\va 1\leq i \leq n',\ (s_{i1}',\ldots,s_{im_i'}')\in S_i,\ t_{i1}'=\ldots = t_{im_i'}'=0 \}.
		\end{align*} 
		We note that there is no relation between $\vec{x},\vec{u}$ and $\vec{y},\vec{v}$, so we get that $S = \pi_1(S)\ti \pi_2(S) = S_1 \ti S_2$. With \th\ref{diophprojectionR1R2} the sets $S_1$ and $S_2$ are Diophantine
	\end{proof}

	\subsection{Diophantine Maps} \phantom{=}\\
	We recall the definition of a Diophantine map from the introduction:
	\begin{defi}
		Let $(R_1,\L_1)$ and $(R_2,\L_2)$ be sets with a language. A map $d: R_1 \ra R_2$ is called a \textbf{Diophantine map} if for all Diophantine sets $S \subset R_1^k$ the set $d(S) \subset R_2^k$ is Diophantine.
	\end{defi}
	
	\begin{prop} \th\label{compdiop=dioph} 
		Let $(R_1,\L_1),(R_2,\L_2)$ and $(R_3,\L_3)$ be sets with a language. Let $d_1: R_1 \ra R_2$ and $d_2: R_2 \ra R_3$ be Diophantine maps, then $d_2 \circ d_1: R_1 \ra R_3$ is Diophantine. \qed 
	\end{prop}
	
	\begin{lemma} \th\label{diophcartesianmapR1R1'R2R2'} Let $d_1:R_1 \ra R_1'$ and $d_2: R_2 \ra R_2'$ be Diophantine maps. Let $\L_i$ respectively $\L_i'$ be the language on $R_i$ respectively $R_i'$ and equip $R_1 \ti R_2$ with $\L_1 \sqcup \L_2$ and $R_1'\ti R_2'$ with $\L_1' \sqcup \L_2'$. Suppose that $\emptyset \subset R_1$ and $\emptyset \subset R_2$ are both a Diophantine set, or both not a Diophantine set. Then $d_1 \ti d_2: R_1 \ti R_2 \ra R_1'\ti R_2'$ defined by $(d_1\ti d_2)(x,y) = (d_1(x),d_2(y))$ is a Diophantine map. \qed 
	\end{lemma}
	
	\begin{lemma} \th\label{diophRkproj} 
		Let $(R,\L)$ be a set with a language. Let $\L^{(k)}$ be the language on $R^{k}$. Then $\pi_i: R^k \ra R$ defined by $\pi_i(\vec{x}) = x_i$ is a Diophantine map. \qed 
	\end{lemma}
	
	\begin{lemma} \th\label{diophR1tR2proj}
		Let $(R_1,\L_1)$ and $(R_2,\L_2)$ be sets with a language such that $\emptyset \subset R_1$ and $\emptyset \subset R_2$ are both a Diophantine set or both not a Diophantine set. Equip $R_1\ti R_2$ with $\L_1 \sqcup \L_2$.\\ 
		Then $\pi_1: R_1\ti R_2 \ra R_1,\ (x,y) \ra x$ and $\pi_2: R_1 \ti R_2 \ra R_2,\ (x,y) \ra y$ are Diophantine maps. \qed 
	\end{lemma} 
	
	\begin{lemma} \th\label{diophmapfinite}
		Let $(R_1,\L_1)$ be a set with a language and let $(R_2,\L_2)$ be an integral domain with $\L_2$ the language of rings with set of constants $S_c(\L_2)$. Suppose that $R_1$ or $R_2$ is finite and $S_c'(\L_2)= R_2$. Then $d: R_1 \ra R_2$ is a Diophantine map.
	\end{lemma}
	
	\begin{proof}
		The proof of this relies on the fact that every finite subset of $R_2$ is Diophantine.
	\end{proof} 
	
	\begin{lemma} \th\label{diophmapf(di)R1rR2} Let $d_1,\ldots,d_r: R_1 \ra R_2$ be Diophantine maps and let $f \in \L_2$ be an $r$-ary function. Let $\L_i$ be the language on $R_i$ and equip $R_1^r$ with $\sqcup_{j=1}^r \L_1$. The map $f(d_1,\ldots,d_r): R_1^r \ra R_2$ defined by $\vec{x} \ra f(d_1(x_1),\ldots,d_r(x_r))$ is Diophantine. \qed 
	\end{lemma}
	
	\begin{remark}
		Let $(R,\L)$ be a set with a language. We can see an $\L$-term $t$ that contains exactly the variables $y_1,\ldots,y_r$ as an $r$-ary function, by $t: R^r \ra R$, $\vec{y} \ra t(\vec{y})$. If we do this, the previous lemma holds with $f$ an $\L_2$-term.
	\end{remark}
	
	\begin{lemma} \th\label{diophmapcases} 
		Let $(R_1,\L_1)$  and $(R_2,\L_2)$ be sets with a language such that $\L_2$ is one-sided, $\cdot \in \L_2$ and $a\cdot b = 0$ implies $a=0$ or $b=0$ for all $a,b \in R_2$. Let $S_1,\ldots,S_n \subset R_1$ be Diophantine sets such that $R_1 = S_1 \sqcup \ldots \sqcup S_n$. If $d_1,\ldots,d_n:R_1 \ra R_2$ are Diophantine maps, then the map $d:R_1 \ra R_2$ defined by $d(x) = d_i(x)$ if $x \in S_i$ is a Diophantine map. \qed 
	\end{lemma}
	
	\begin{remark} \th\label{diophmapfromsubset} 
		Let $S \subset (R_1,\L_1)$ be a Diophantine set such that $S^c = R_1\bs S$ is also Diophantine. Then the previous lemma makes it possible to speak of Diophantine maps from $S \ra (R_2,\L_2)$ since we can extend any Diophantine map $d:S \ra (R_2,\L_2)$ in a Diophantine way to $R_1 \ra R_2$, by defining $d(x) = 0$ if $x \notin S$. 
	\end{remark} 
	
	\begin{thm}[\th\ref{introDiophMap}] \th\label{diophset->diophset} \th\label{diophinjdef} \th\label{diophinjdefLwithR} 
		Let $(R_1,\L_1)$ and $(R_2,\L_2)$ be sets with a language. Let $d: R_1 \ra R_2$ be an injective map such that:
		\begin{itemize}
			\item For all constants $c \in \L_1$, $d(c) \in S_c'(\L_2)$, i.e. $d(c)$ can be defined by $\L_2$;
			\item For all functions $f: R_1^k \ra R_1$ of $\L_1$ the set $\{(d(x_1),\ldots,d(x_k),d(f(x_1,\ldots,x_k))) \mid \vec{x} \in R_1^k\}$ is Diophantine;
			\item For all relations $S \subset R_1^k$ of $\L_1$ the set $d(S)$ is Diophantine;
			\item The set $d(R_1) \subset R_2$ is Diophantine.
		\end{itemize}
		Then $d$ is a Diophantine map.
	\end{thm}
	
	\begin{proof}
		We will first prove by induction with respect to the number of functions in an $\L_1$-term that for all $\L_1$-terms $t$ the set $\{(d(\vec{x}),d(t(\vec{x})))\mid x \in R_1^k\})$ is Diophantine.\\
		\textit{Induction basis} If $t$ contains no functions, then $t$ is a variable or a constant. In the first case we get $\{(d(x),d(x))\mid x \in R_1\}) = \{(u,v) \in R_2^2\mid u,v \in d(R_1),\ u = v\}$, which is a Diophantine set, since by assumption $d(R_1)$ is a Diophantine set. In the second case we get $\{(d(\vec{x}),d(c))\mid \vec{x} \in R_1^0 \} = \{d(c)\}$, which is a Diophantine set per assumption.\\
		\textit{Induction hypothesis} Suppose that for all $\L_2$-terms $t$ where $t$ contains $n$ or less functions the set $\{(d(\vec{x}),d(t(\vec{x})))\mid x \in R_1^k\}$ is a Diophantine set. \\
		\textit{Induction step} Let $t$ be a term with $n+1$ functions. Then $t$ is of the form $t= f(t_1,\ldots,t_m)$ with $f \in \L_1$ and $t_i$ terms of $\L_1$ with $n$ or less functions. With the induction hypothesis we get that the sets $\{(d(\vec{x_i}),d(t_i(\vec{x_i})))\mid \vec{x} \in R_1^{k_i}\})$ are Diophantine sets. Since the Cartesian product of Diophantine sets is again Diophantine by \th\ref{diophcartRR} we get that also the sets $\{(d(\vec{x}),d(t_i(\vec{x})))\mid \vec{x} \in R_1^k\})$, so with $k$ instead of $k_i$ variables, are Diophantine. Then since $d$ is injective we get
		\begin{align*}
			\{(d(\vec{x}),d(t(\vec{x})))\mid \vec{x}\in R_1^k\} = \{(d(\vec{x}),d(f(t_1(\vec{x}),\ldots,t_m(\vec{x}))))\mid \vec{x} \in R_1^k\} \\
			= \{(\vec{u},v)\in R_2^{k+1}\mid \ex (\vec{u_i},v_i)\in \{(d(\vec{x}),d(t_i(\vec{x})))\mid \vec{x} \in R_1^k\},\ \va 1\leq i \leq m,\ \vec{u} = \vec{u_i},\\ 
			(v_1,\ldots,v_n,v) \in \{(d(\vec{y}),d(f(\vec{y})))\mid \vec{y} \in R_1^k\}\}.
		\end{align*}
		Since by assumption also $\{(d(\vec{y}),d(f(\vec{y})))\mid \vec{y} \in R_1^k\}$ is a Diophantine set, this implies that the set $\{(d(\vec{x}),d(t(\vec{x})))\mid \vec{x}\in R_1^k\})$ is Diophantine, which proves the induction hypothesis. \\
		Suppose now that $S \subset R_1^k$ is a Diophantine set. Then we can write
		$$S = \{\vec{x} \in R_1^k\mid \ex \vec{y} \in R_1^{l},\ \va 1\leq i \leq n,\ (t_{i1}(\vec{x},\vec{y}),\ldots,t_{im_i}(\vec{x},\vec{y})) \in S_i\}$$
		where $t_{ij}$ are $\L_1$-terms and $S_i$ are relations of $\L_1$, which can include equality. We get
		\begin{align*}
			d(S) &= \{\vec{u} \in R_2^k\mid \ex \vec{v} \in R_2^l,\ \ex z_{ij} \in R_2,\\
			&\phantom{===} \vec{u} \in d(R_1^k),\ \vec{v} \in d(R_1^l),\ z_{ij} \in d(R_1) \\
			&\phantom{====} \va 1 \leq i \leq n,\ \va 1 \leq j \leq m_i,\ (\vec{u},\vec{v},z_{ij}) \in \{(d(\vec{x}),d(\vec{y}),d(t_{ij}(\vec{x},\vec{y})))\} \\
			&\phantom{=====} \va 1\leq i \leq n,\ \vec{z_i} \in d(S_i) \}. 
		\end{align*} 
		Since $d$ is an injective map, the second line forces $\vec{u} = d(\vec{x})$ and $\vec{v} = d(\vec{y})$ for unique $\vec{x} \in R_1^k$ and $\vec{y}\in R_2^l$. The third line forces $z_{ij} = d(t_{ij}(\vec{x},\vec{y}))$. If we combine this with the last line we get that $(d(t_{i1}(\vec{x},\vec{y})),\ldots,d(t_{im_i}(\vec{x},\vec{y}))) \in d(S_i)$. Since $d$ is still injective, this gives $(t_{i1}(\vec{x},\vec{y}),\ldots,t_{im_i}(\vec{x},\vec{y})) \in S_i$.\\
		Since we have proven that $\{(d(\vec{x}),d(\vec{y}),d(t_{ij}(\vec{x},\vec{y})))\mid \vec{x} \in R_1^k,\ \vec{y} \in R_2^l\}$ are Diophantine sets and $d(R_1)$ and the $d(S_i)$ are Diophantine sets by assumption, we get that $d(S)$ is a Diophantine set. This implies that $d$ is a Diophantine map.
	\end{proof} 
	
	\begin{remark}
		The condition in \th\ref{diophinjdef} that $d(R_1)$ is a Diophantine set gives that 
		$$d(\{(x,y) \in R_1\mid x=y\}) = \{(x,y) \in R_2^2\mid x \in d(R_1),\ x=y\}$$ is a Diophantine set. This means that we can omit to check the condition for the relation equality in \th\ref{diophinjdef}. 
	\end{remark} 
	
	\begin{example} \th\label{diophexZZx->x+n}
		With \th\ref{diophinjdef}, we can prove that the injective map $d: \Z \ra \Z$, $x \ra x+n$ for some fixed $n\in \Z$ is a Diophantine map.\\
		We have $d(\Z) = \Z = \{x \in \Z\}$ is Diophantine. We have $d(0) = n$ and $d(1) = n+1$, which are in $S_c'(\Z) = \Z$. Further, we get
		\begin{align*} 
			\{(d(x),d(y),d(x+y))\mid x,y \in \Z\} &= \{(x+n,y+n,x+y+n)\mid x,y \in \Z\}\\ 
			&= \{(x',y',z')\in \Z^3 \mid x'+y'-z'-n=0\},
		\end{align*} 
		which is also Diophantine, as $n \in S_c'(\Z)$. Finally, we have
		\begin{align*} 
			&\phantom{=} \{(d(x),d(y),d(x\cdot y))\mid x,y \in \Z\} = \{(x+n,y+n,x\cdot y+n)\mid x,y \in \Z\} \\
			&=\{(x',y',z') \in \Z^3 \mid x'\cdot y'-n\cdot x'-n\cdot y'+ n^2+n-z'=0\}.
		\end{align*} 
		This is also a Diophantine set, so $d$ satisfies the condition of \th\ref{diophset->diophset}, so it is Diophantine.
	\end{example}

	\begin{example} \th\label{diophinclNZ}
		Equip $\N$ and $\Z$ with the language of rings. Then the inclusion $\i: \N \ra \Z$ is a Diophantine map.
	\end{example} 
	
	\begin{proof}
		We have that $\i$ is injective, so we are going to use \th\ref{diophset->diophset}. \\
		By Lagrange's four squares theorem, we get that 
		$$\i(\N) = \N = \{x \in \Z \mid \ex a,b,c,d \in \Z,\ a^2+b^2+c^2+d^2-x =0\}$$ is a Diophantine set. \\
		We have that $\i(0) = 0 \in S_c(\Z)$ and $\i(1)=1 \in S_c(\Z)$. \\
		For $* \in\{+,\cdot\}$ we get that 
		$$\{(\i(a),\i(b),\i(a*b))\mid a,b \in \N\} = \{(x,y,z) \in \Z^3\mid x,y \in \N,\ z = x*y\}.$$
		This is a Diophantine subset of $\Z^3$, since $\N \subset \Z$ is a Diophantine set.\\
		These are all conditions of the theorem, which implies that $\i$ is a Diophantine map.
	\end{proof} 
	
	\begin{remark}
		If $d$ is not injective, the statement of \th\ref{diophset->diophset} is false. Equip $\Z$ with the language $\L = \{+,0,1\}$ and $\Z^2$ with $\L \sqcup \L$. We will show that $$d: \Z^2 \ra \Z,\ (a,b) \ra a^2b$$ meets all requirements of the theorem except the injectivity. Denote by $+_i$ and $-_i$ the elements of $\L\sqcup \L$ that are addition respectively subtraction in the $i$-th coordinate. Then we show that the image of the Diophantine set 
		$$S = \{\vec{a} \in \Z^2\mid a -_2 1 = 0\} = \{(a_1,1) \in \Z^2\}$$ is not a Diophantine set. This implies that $d$ is not Diophantine.\\
		We have that $S_c'(\L) = \Z$, so $\{d(c)\} = \{x \in \Z\mid x-d(c)=0\}$ is Diophantine for all constants $c \in \L\sqcup \L$. Next, we have $d(\Z^2) = \Z$, which is a Diophantine set. We have
		\begin{align*} 
			&\phantom{=} \{(d(\vec{x}),d(\vec{y}),d(\vec{x}+_1\vec{y}))\mid \vec{x},\vec{y} \in \Z^2\} 
			= \{(d(\vec{x}),d(\vec{y}),d((x_1+y_1,0)))\mid \vec{x},\vec{y} \in \Z^2\} \\
			&= \{(x_1^2x_2,y_1^2y_2, 0)\mid x_1,x_2,y_1,y_2 \in \Z\}
			= \Z^2\ti \{0\} = \{\vec{x} \in \Z^3\mid x_3=0\};\\
			&\phantom{=} \{(d(\vec{x}),d(\vec{y}),d(\vec{x}+_2\vec{y}))\mid \vec{x},\vec{y} \in \Z^2\} 
			= \{(d(\vec{x}),d(\vec{y}),d((0,x_2+y_2)))\mid \vec{x},\vec{y} \in \Z^2\} \\
			&= \{(x_1^2x_2,y_1^2y_2, 0)\mid x_1,x_2,y_1,y_2 \in \Z\}
			= \Z^2\ti \{0\} = \{\vec{x} \in \Z^3\mid x_3=0\}.
		\end{align*} 
		These are Diophantine sets, so $d$ fulfils all conditions of \th\ref{diophset->diophset} except the injectivity.\\
		Since we only have the operation $+$ in $\L$, an arbitrary Diophantine subset $S'$ of $(\Z,\L)$ is given by
		$$S' = \left\{x \in \Z \mid \ex y_1,\ldots,y_n\text{ s.t. }\va\ 1\leq j \leq k\ a_{0j}x + \left(\sum_{i=1}^n a_{ij}y_i \right) -b_j=0\right\},$$
		for some fixed $a_{ij},b_j \in \Z$. This set of equations is equivalent to the matrix equation $A\cdot (x,\vec{y})^T = \vec{b}^T$, where $A = (a_{ij})_{i=0,j=1}^{n,k}$. From linear algebra over the integers, see \cite{LinAlgZ}, we know that an integer solution to this system, if an integer solution exists, is given by $(x,\vec{y})^T = V\cdot (c_0,\ldots,c_k,z_{k+1},\ldots,z_n)^T$ for some integer matrix $V$, constants $c_i \in \Z$ and arbitrary integers $z_i$. So $S'$ is of the form $c_0'+c_1'\Z + \ldots + c_{k'}'\Z$ for some $k'\in \Z_{\geq 0}$, $c_i'\in \Z$ constants, or $S'=\emptyset$. \\
		We have that $d(S) = d(\{(a_1,1)\in \Z^2\}) = \{a_1^2\mid a_1 \in \Z\}$ is not of this form, so $d(S)$ is not Diophantine.
	\end{remark}
	
	\subsection{Maps Given by Diophantine Expressions}
	
	\begin{defi}[\th\ref{introDefGivenbyDiophExpr}] Let $(R_1,\L_1)$ and $(R_2,\L_2)$ be sets with a language and equip $R_1^k$ with $\L_1^{(k)}$ and $R_2^l$ with $\L_2^{(l)}$. Let $d_1: R_1 \ra R_2$ and $d_2: R_1^k \ra R_2^l$ be maps. We say that $d_2$ is \textbf{given by a Diophantine expression} in $d_1$ if there exist basic $\L_2$-formulas $f_1,\ldots,f_n$ such that 
		$$d_2(\vec{x}) = \vec{z} \ \Leftrightarrow\ \ex \vec{y} \in R^{l'},\ \va 1\leq i \leq n,\ f_i(d_1(\vec{x}),\vec{y},\vec{z})$$ 
		If $k=l=1$, we will denote this relation with $d_2 \geq_D d_1$ or $d_1 \leq_D d_2$.
	\end{defi}
	
	\begin{lemma} \th\label{dioph<Dreflexive}
		Let $d:R_1 \ra R_2$ be a map. Then $d \geq_D d$, so $\geq_D$ is reflexive. \qed 
	\end{lemma} 
	
	\begin{lemma} \th\label{dioph<Dtransitive} 
		Let $d_1,d_2,d_3:R_1 \ra R_2$ be maps such that $d_1 \leq_D d_2$ and $d_2 \leq_D d_3$ and $d_2(R_1)$ is a Diophantine set. Then $d_1 \leq_D d_3$, so $\leq_D$ is transitive. \qed 
	\end{lemma}
	
	\begin{prop} \th\label{diophd1dioph+d1<Dd2=>d2dioph} 
		Let $(R_1,\L_1)$ and $(R_2,\L_2)$ be sets with a language. Equip $R_1^k$ and $R_2^{k'}$ with $\L_1^{(k)}$ respectively $\L_2^{(k')}$. Let $d_1:R_1 \ra R_2$ and $d_2: R_1^k\ra R_2^{k'}$ be maps. If $d_1$ is an injective Diophantine map and $d_2$ is given by a Diophantine expression in $d_1$, then $d_2$ is a Diophantine map.
	\end{prop} 
	
	\begin{proof}
		We have that $d_2$ is given by a Diophantine expression in $d_1$, so there exist basic $\L_2$-formulas $f_1,\ldots, f_n$ such that
		$$d_2(\vec{x})=\vec{z} \LRa \ex \vec{w} \in R_2^l ,\ \va 1\leq i \leq n,\ f_i(d_1(\vec{x}),\vec{w},\vec{z}).$$
		Let $S \subset R_1^{kk'}$ be a Diophantine set. Then there exist $\L_1^{(k)}$-terms $t_{ij}$ and $\L_1^{(k)}$-relations $S_i$ such that 
		$$S = \{\vec{x} \in R_1^{kk'}\mid \ex \vec{y} \in R_1^{kl'},\ \va 1\leq i \leq n',\ (t_{i1}(\vec{x},\vec{y}),\ldots,t_{im_i}(\vec{x},\vec{y})) \in S_i\}.$$
		We have that the set $\{(\vec{x},\vec{y})\mid (t_{i1}(\vec{x},\vec{y}),\ldots,t_{im_i}(\vec{x},\vec{y})) \in S_i\}$ is a Diophantine set. Since $d_1$ is a Diophantine map, this gives that $$T_i = \{(d_1(\vec{x}),d_1(\vec{y}))\mid (t_{i1}(\vec{x},\vec{y}),\ldots,t_{im_i}(\vec{x},\vec{y})) \in S_i\}$$
		is a Diophantine set, so there exist $\L_2$-formulas such that
		$$T_i = \{(\vec{u},\vec{v})\mid \ex \vec{y_i} \in R_2^{ll_i},\ \va 1\leq j \leq k_i,\ f_{ij}(\vec{u},\vec{v},\vec{y_i}) \}.$$
		Since $d_1$ is injective, this gives that
		$$(t_{i1}(\vec{x},\vec{y}),\ldots,t_{im_i}(\vec{x},\vec{y})) \in S_i \LRa \ex \vec{y_i} \in R_2^{ll_i},\ \va 1\leq j \leq k_i,\ f_{ij}(d_1(\vec{x}),d_1(\vec{y}),\vec{y_i}) $$
		Then we get that
		\begin{align*} 
			d_2(S) &= \{d_2(\vec{x})\mid \vec{x}\in R_1^{kk'},\ex \vec{y}\in R_1^{kl'},\ \va 1\leq i\leq n',\ (t_{i1}(\vec{x},\vec{y}),\ldots,t_{im_i}(\vec{x},\vec{y})) \in S_i\}\\
			&= \{\vec{z} \in R_2^k\mid \ex \vec{u} \in d_1(R_1^k)^{k'},\ \ex \vec{v} \in d_1(R_1^k)^{l'},\\ 
			&\phantom{===} \va 1\leq j\leq k,\ \ex \vec{w_j} \in R_2^{kl'},\ \va 1\leq i \leq n,\ f_i(\vec{u_j},\vec{w_j},\vec{z_j}),\\
			&\phantom{===} \va 1\leq i \leq n',\ \ex \vec{y_i} \in R_2^{ll_i},\ \va 1\leq j \leq k_i,\ f_{ij}(\vec{u},\vec{v},\vec{y_i})\}.
		\end{align*} 
		This holds since $\vec{u} \in d_1(R_1^k)^{k'}$ and $\vec{v} \in d_1(R_1^k)^{l'}$ imply that $u = d_1(\vec{x})$ and $v = d_1(\vec{y})$ for unique $\vec{x}\in R_1^{kk'}$, $\vec{y} \in R_1^{kl'}$, as $d_1$ is injective. Then the formulas $f_i$ force $\vec{z_j} = d(\vec{x_j})$ and the formulas $f_{ij}$ force the vectors $\vec{x}$ and $\vec{y}$ to satisfy the defining formulas of $S$, so $\vec{x} \in S$.\\
		Since $d_1$ is Diophantine, the set $d_1(R_1^k)$ is Diophantine, so $d_2(S)$ is a Diophantine set. This implies that $d_2$ is a Diophantine map.
	\end{proof}
	
	\begin{cor} \th\label{diophid<=Dd=>ddioph}
		Let $(R,\L)$ be a set with a language. Equip $R^k$ and $R^l$ with $\L^{(k)}$ respectively $\L^{(l)}$. If $d: R^k \ra R^l$ is such that there exist basic formulas $f_1,\ldots,f_n$ such that 
		$$d(\vec{x})= \vec{z} \LRa \ex \vec{y} \in R^{l'},\ \va 1\leq i \leq n,\  f_i(\vec{x},\vec{y},\vec{z}),$$ 
		then $d$ is Diophantine. In other words, if $d$ is given by a Diophantine expression, it is a Diophantine map.
	\end{cor}
	
	\begin{lemma} \th\label{diophd2dioph+d1<Dd2=>d1dioph}
		Let $d_1,d_2: R_1 \ra R_2$ be injective maps such that $d_1 \leq_D d_2$. If $d_1(R_1)$ is a Diophantine set and $d_2$ is a Diophantine map, then $d_1$ is a Diophantine map.
	\end{lemma}
	
	\begin{proof}
		As $d_1 \leq_D d_2$, there exist basic $\L_2$-formulas $f_1,\ldots,f_n$ such that
		$$d_2(x)=z \LRa \ex y \in R_2^l,\ \va 1\leq i \leq n,\ f_i(d_1(x),\vec{y},z).$$
		Since $d_1$ is injective we can substitute $x = d_1^{-1}(u)$, which gives
		$$d_2(d_1^{-1}(u)) = z \LRa u \in d_1(R_1),\ \ex y\in R_2^l,\ \va 1\leq i \leq n,\ f_i(u,\vec{y},z).$$
		Let $S \subset R_1^k$ be a Diophantine set. Then as $d_1$ and $d_2$ are injective, we get
		\begin{align*}
			d_1(S) &= \{d_1(\vec{x})\mid \vec{x} \in S\}\\
			&= \{d_1(\vec{x})\mid d_2(\vec{x}) \in d_2(S)\}\\
			&= \{\vec{u}\in R_2^k \mid \vec{u} \in (d_1(R_1))^k,\ \ex \vec{z} \in d_2(S),\  d_2(d_1^{-1}(\vec{u})) = \vec{z} \} \\
			&= \{\vec{u} \in R_2^k\mid \vec{u} \in (d_1(R_1))^k,\ \ex \vec{z} \in d_2(S),\ \va 1\leq j \leq k,\ \ex \vec{y_j},\ \va 1\leq i \leq n,\ f_i(u_j,\vec{y_j},z_j)\}. 
		\end{align*}
		The last set is Diophantine, as $d_2$ is a Diophantine map and we require $d_1(R_1)$ to be a Diophantine set. So $d_1$ is a Diophantine map.
	\end{proof}
	
	\begin{lemma}
		Let $(R_1,\L_1)$ and $(R_2,\L_2)$ be sets with a language and let $d:R_1 \ra R_2$ be a map. Suppose that $d_1,\ldots,d_r: R_1 \ra R_2$ are given by Diophantine expressions in $d$. For $f \in \L_2$ the map $f(\vec{d}):R_1 \ra R_2,\ f(\vec{d})(x)= f(d_1(x),\ldots,d_r(x))$ is given by a Diophantine expression in $d$.\\
		If moreover, $d$ is an injective Diophantine map, then the $d_i$ and $f(\vec{d})$ are Diophantine maps. \qed 
	\end{lemma}
	
	\begin{remark}
		If $\L_2$ contains the language of rings and $R_2$ is a ring, we get for a fixed $d: R_1 \ra R_2$ a ring-structure on the set $\{d':R_1 \ra R_2\mid d\leq_D d'\}$. If $d$ is injective and Diophantine, this is a subset of the set of Diophantine maps from $R_1$ to $R_2$. This leads to the question of when these sets are equal, as we then have a ring structure on our set of Diophantine maps. We will consider this question in Paragraph \ref{parRecvsDioph}. 
	\end{remark}
	
	\begin{lemma} \th\label{diophcartmapR12R22}
		Let $(R_1,\L_1)$ and $(R_2,\L_2)$ be sets with a language. Let $d,d_1,d_2: R_1 \ra R_2$ be maps such that $d \leq_D d_1$ and $d\leq_D d_2$. Then the map $d_1\ti d_2: (R_1^2,\L_1^{(2)}) \ra (R_2^2,\L_2^{(2)}),\ (x,y) \ra (d_1(x),d_2(y))$ is given by a Diophantine expression in $d$.\\ 
		If moreover, $d$ is an injective Diophantine map, then $d_1\ti d_2$ is a Diophantine map. \qed 
	\end{lemma} 
	
	\begin{cor} \th\label{diophcartR1R22} 
		Suppose that the conditions of \th\ref{diophcartmapR12R22} hold. Then the map $(d_1,d_2):R_1 \ra R_2^2$ defined by $x \ra (d_1(x),d_2(x))$ is given by a Diophantine expression in $d$. If furthermore $d$ is an injective Diophantine map, then $(d_1,d_2)$ is a Diophantine map. \qed 
	\end{cor}
	
	\begin{remark}
		The condition that $d_1$ and $d_2$ are given by Diophantine expressions in some injective Diophantine map is necessary in \th\ref{diophcartR1R22}. We will construct a counter-example where $d_1$ and $d_2$ are Diophantine, but $d_1 \ti d_2$ is not. \\
		Let $R = \Z^2$ and let $\L = \{0,1,+,\cdot\}$ with the coordinate-wise interpretation. Let $d_1: R \ra R$ be the identity and let $d_2:R \ra R, (x,y) \ra (y,x)$. As $(R,\L)$ is symmetric, both $d_1$ and $d_2$ are Diophantine maps. We take the Diophantine set 
		$$S = \{(\vec{x},\vec{y}) \in R^2\mid \vec{x} -\vec{y} = 0\} = \{((a,b),(a,b))\mid a,b \in \Z\}.$$
		Then $S'= (d_1\ti d_2)(S) = \{((a,b),(b,a))\mid a,b \in \Z\}$. Suppose that $S'$ is Diophantine. Then there exist polynomials $P_1,\ldots,P_n$ such that 
		\begin{align*}
			S'&= \{(\vec{x},\vec{y}) \in (\Z^2)^2\mid \ex ((a_1,b_1),\ldots,(a_k,b_k)) \in (\Z^2)^k,\ \va 1\leq i\leq n,\  P_i(\vec{x},\vec{y},(a_1,b_1),\ldots,(a_k,b_k)) = 0\}\\
			&= \{(\vec{x},\vec{y}) \in (\Z^2)^2\mid \ex \vec{a},\vec{b} \in \Z^k,\ \va 1\leq i\leq n,\  P_i(x_1,y_1,\vec{a}) = P_i(x_2,y_2,\vec{b}) = 0\}.
		\end{align*}  
		As $((0,0),(0,0))$ and $((1,1),(1,1))$ are in $S'$, there exist $\vec{a},\vec{b} \in \Z^k$ and $\vec{c},\vec{d}\in \Z^k$ such that for all $1\leq i \leq n$ we have $P_i(0,0,\vec{a}) = P_i(0,0,\vec{b}) = P_i(1,1,\vec{c}) = P_i(1,1,\vec{d}) = 0$. But then also $((0,1),(0,1))\in S'$ as for all $1\leq i \leq n$ we have $P_i(0,0,\vec{a}) = P_i(1,1,\vec{c}) = 0$. This gives a contradiction, so $S'= (d_1\ti d_2)(S)$ is not a Diophantine set and $d_1 \ti d_2$ is not a Diophantine map.
	\end{remark}

	\subsection{Diophantine Equivalence Maps}

	\begin{defi} [\th\ref{introDefDiophEquiv}]
		Let $(R_1,\L_1)$ and $(R_2,\L_2)$ be sets with a language. Let $\sim$ be an equivalence relation on $R_2$ and let $\pi: R_2 \ra R_2\modsim$ be the projection. If $d: R_1 \ra R_2\modsim$ is a map such that for all Diophantine sets $S \subset R_1^k$ the set 
		$$\ol{d}(S) = \{\vec{x} \in R_2^k \mid \va 1\leq i \leq k,\ \pi(x_i) \in d(S)\}$$
		is Diophantine, then $d$ is a \textbf{Diophantine equivalence map}. 
	\end{defi}
	
	\begin{thm}[\th\ref{introDiophEquiv}] \th\label{H10Pdinjdiophequiv} 
		Let $(R_1,\L_1)$ and $(R_2,\L_2)$ be sets with a language. Let $\sim$ be an equivalence relation on $R_2$ and let $\pi:R_2 \ra R_2\modsim$ be the projection. Let $d: R_1 \ra R_2\modsim$ be an injective map such that
		\begin{itemize}
			\item for all constants $c \in \L$, $\ol{d}(c) \subset R_2$ is a Diophantine set;
			\item for all $k$-ary functions $f \in \L$, $\ol{d}(\{(\vec{x},f(\vec{x}))\mid \vec{x} \in R_1^k\})$ is a Diophantine set;
			\item for all relations $S \subset \L$, $\ol{d}(S)$ is a Diophantine set;
			\item the set $\ol{d}(R_1)$ is Diophantine;
			\item the set $\{(x,y) \in R_2^2\mid x \sim  y\}$ or the set $\{(x,y)\in R_2^2 \mid x,y \in \ol{d}(R_1),\ x\sim y\}$ is Diophantine.
		\end{itemize}
		Then $d$ is a Diophantine equivalence map.
	\end{thm} 
	
	\begin{proof}
		We will first prove by induction on the number of functions in an $\L_1$-term that for all $\L_1$-terms $t$ the set $\ol{d}(\{(\vec{x},t(\vec{x}))\mid x \in R_1^k\})$ is Diophantine.\\
		\textit{Induction basis} If $t$ contains no functions, then $t$ is a variable or a constant. In the first case we get $\ol{d}(\{(x,x)\mid x \in R_1\}) = \{(u,v) \in R_2^2\mid u,v \in \ol{d}(R_1),\ u\sim v\}$. By the last condition or the last two conditions and \th\ref{diophintersection} this is a Diophantine set. In the second case we get $\ol{d}(\{(\vec{x},c) \in R_1^0\ti R_1 \}) = \ol{d}(c)$, which is a Diophantine set per assumption.\\
		\textit{Induction hypothesis} Suppose that for all $\L_2$-terms $t$ where $t$ contains $n$ or less functions the set $\ol{d}(\{(\vec{x},t(\vec{x}))\mid x \in R_1^k)\})$ is a Diophantine set. \\
		\textit{Induction step} Let $t$ be a term with $n+1$ functions. Then $t$ is of the form $t= f(t_1,\ldots,t_m)$ with $f \in \L_1$ and $t_i$ terms of $\L_1$ with $n$ or less functions. By the induction hypothesis we get that the sets $\ol{d}(\{(\vec{x_i},t_i(\vec{x_i}))\mid \vec{x} \in R_1^{k_i}\})$ are Diophantine. Since the Cartesian product of Diophantine sets is again Diophantine by \th\ref{diophcartRR} we get that also the sets $\ol{d}(\{(\vec{x},t_i(\vec{x}))\mid \vec{x} \in R_1^k\})$, so with $k$ instead of $k_i$ variables, are Diophantine. Then
		\begin{align*}
			\ol{d}(\{(\vec{x},t(\vec{x}))\mid \vec{x}\in R_1^k\}) = \ol{d}(\{(\vec{x},f(t_1(\vec{x}),\ldots,t_m(\vec{x})))\mid \vec{x} \in R_1^k\}) \\
			= \{(\vec{u},v)\in R_2^{k+1}\mid \ex (\vec{u_i},v_i)\in \ol{d}(\{(\vec{x},t_i(\vec{x}))\mid \vec{x} \in R_1^k\}),\ \va 1\leq i \leq m,\ \vec{u} = \vec{u_i},\\ 
			(v_1,\ldots,v_n,v) \in \ol{d}(\{(\vec{y},f(\vec{y}))\mid \vec{y} \in R_1^k\})\}.
		\end{align*}
		Since by assumption also $\ol{d}(\{(\vec{y},f(\vec{y}))\mid \vec{y} \in R_1^k\})$ is a Diophantine set, we get that $\ol{d}(\{(\vec{x},t(\vec{x}))\mid \vec{x}\in R_1^k\})$ is a Diophantine set, which proves the induction hypothesis. \\
		Suppose now that $S \subset R_1^k$ is a Diophantine set. Then we can write
		$$S = \{\vec{x} \in R_1^k\mid \ex \vec{y} \in R_1^{l},\ \va 1\leq i \leq n,\ (t_{i1}(\vec{x},\vec{y}),\ldots,t_{im_i}(\vec{x},\vec{y})) \in S_i\}$$
		where $t_{ij}$ are $\L_1$-terms and $S_i$ are relations of $\L_1$, which can include equality. We get
		\begin{align*}
			\ol{d}(S) &= \{\vec{u} \in R_2^k\mid \ex \vec{v} \in R_2^l,\ \ex z_{ij} \in R_2,\\
			&\phantom{===} \vec{u} \in \ol{d}(R_1^k),\ \vec{v} \in \ol{d}(R_1^l),\ z_{ij} \in \ol{d}(R_1) \\
			&\phantom{====} \va 1 \leq i \leq n,\ \va 1 \leq j \leq m_i,\ (\vec{u},\vec{v},z_{ij}) \in \ol{d}(\{(\vec{x},\vec{y},t_{ij}(\vec{x},\vec{y}))\}) \\
			&\phantom{=====} \va 1\leq i \leq n,\ \vec{z_i} \in \ol{d}(S_i) \}. 
		\end{align*} 
		Let $[\vec{u}]$ denote the equivalence class of $\vec{u}$ under $\sim$. Since $d$ is an injective map, the second line forces $[\vec{u}] = [d(\vec{x})]$ and $[\vec{v}] = [d(\vec{y})]$ for unique $\vec{x} \in R_1^k$ and $\vec{y}\in R_2^l$. The third line forces $[z_{ij}] = [d(t_{ij}(\vec{x},\vec{y}))]$. Combining this with the last line gives $([d(t_{i1}(\vec{x},\vec{y}))],\ldots,[d(t_{im_i}(\vec{x},\vec{y}))]) \in \ol{d}(S_i)$. Since $d$ is still injective, this gives $(t_{i1}(\vec{x},\vec{y}),\ldots,t_{im_i}(\vec{x},\vec{y})) \in S_i$.\\
		Since we have proven that $\ol{d}(\{(\vec{x},\vec{y},t_{ij}(\vec{x},\vec{y}))\mid \vec{x} \in R_1^k,\ \vec{y} \in R_2^l\}$ are Diophantine sets and $\ol{d}(R_1)$ and $\ol{d}(S_i)$ are Diophantine sets by assumption, we get that $\ol{d}(S)$ is a Diophantine set. This implies that $d$ is a Diophantine equivalence map.
	\end{proof} 
	
	\subsection{Diophantine Models}\phantom{=}\\
	In this paragraph, we will define what a Diophantine model is. Then we will define when two models are similar or equivalent and look at some examples. 
	
	\begin{defi}
		Let $(R_1,\L_1)$ and $(R_2,\L_2)$ be sets with a language. A \textbf{Diophantine model} of $R_1$ in $R_2$ is an injective Diophantine map $d:R_1 \ra R_2$ or an equivalence relation $\sim$ on $R_2$ and an injective Diophantine equivalence map $d: R_1 \ra R_2\modsim$. 
	\end{defi}
	
	\begin{defi}
		Let $(R,\L)$ be a set with a language. A \textbf{Diophantine automorphism} of $R$ is a bijective map $f:R \ra R$ such that $f$ and $f^{-1}$ are Diophantine maps. 
	\end{defi} 
	
	We will first do the case that $d_1$ and $d_2$ do not involve an equivalence relation.
	
	\begin{defi}
		Let $d_1:R_1 \ra R_2$ and $d_2:R_1 \ra R_2$ be models of $(R_1,\L_1)$ in $(R_2,\L_2)$. We say that $d_1$ and $d_2$ are \textbf{similar} if there exists an Diophantine automorphism $f$ of $R_2$ such that $d_1(R_1) = f(d_2(R_1))$. \\
		The models are \textbf{equivalent} if there exist Diophantine automorphisms $f_1$ of $R_1$ and $f_2$ of $R_2$ such that $d_1 = f_2 \circ d_2 \circ f_1$.
	\end{defi} 
	
	\begin{remark}
		Note that the relations of being similar or equivalent models are both equivalence relations. We also have that if $d_1$ and $d_2$ are equivalent models, they are similar since an automorphism is bijective.
	\end{remark} 
	
	\begin{example}
		We can model $\Z$ inside itself with the map $d_n: \Z \ra \Z$ defined by $d_n(x) = x+n$. With \th\ref{diophexZZx->x+n} this map is Diophantine. We have that $d_n(\Z) = \Z$ for all $n$, so it immediately follows that the $d_n$ are similar models. Each $d_n$ is a Diophantine automorphism since its inverse is $d_{-n}$. This implies that $d_n$ and $d_m$ are equivalent since $d_n = d_{n-m}\circ d_m$. \\
		Define for $k\not= 0$ the map $d_k': \Z \ra \Z$ by $d_k'(x) = kx$. Then $d_k'$ is given by a Diophantine expression in the identity, so it is a Diophantine map by \th\ref{diophd1dioph+d1<Dd2=>d2dioph}. \\
		We have that $d_n$ is not similar, so also not equivalent, to $d_k'$ for $k\not= \pm 1$ since we can not bijectively map $\Z \bs d_k'(\Z) = \{x \in \Z \mid k\nmid x\}$ to $\Z \bs d_n(\Z) = \emptyset$. For the same reason $d_1'$ and $d_{-1}'$ are not similar to $d_k'$ for $k\not= \pm 1$. Since $d_{-1}'\circ d_{-1}' = d_1' = id$, the map $d_{-1}'$ is a Diophantine automorphism and $d_1'$ and $d_{-1}'$ are equivalent. We can also show that $d_k'$ and $d_l'$ with $k,l\not= \pm 1$ are equivalent. We have that the sets $\{x \in \Z\mid k\mid x\}$ and $\{x \in \Z \mid k\nmid x\} = \{x \in \Z \mid k\mid \prod_{i=1}^{k-1} x+i\}$ are both Diophantine. Since $\Z$ is a domain and the language of rings is one-sided on $\Z$, we can do definition by cases. \\
		Define $f_1: \{x \in \Z \mid k\mid x\} \ra \{x \in \Z \mid l \mid x\}$ by $f_1(x) = \frac{lx}{k}$. We have that $f_1$ is given by a Diophantine expression in the identity map, namely $f_1(x) = y \LRa \ex q,\ kq = x,\ y = lq$. This implies with \th\ref{diophd1dioph+d1<Dd2=>d2dioph} that $f_1$ is a Diophantine map. \\
		We define $f_2: \{x \in \Z \mid k \nmid x\} \ra \Z$ by writing $x = kq+r$ with $1 \leq r \leq k-1$ and defining $f_2(x) = (k-1)q+r-1$. With Lagrange's four-square theorem, we get that the relation $\leq$ is Diophantine, since $\{(x,y) \mid x\leq y\} = \{(x,y)\mid \ex a,b,c,d \in \Z,\ y-x = a^2 + b^2+c^2+d^2\}$. This implies that $f_2$ is given by a Diophantine expression in the identity map, namely 
		$$f_2(x) = y \LRa \ex q,r,\ kq+r = x,\ 1\leq r \leq k-1,\ y = (k-1)q+r-1.$$ This implies that $f_2$ is a Diophantine map. We have that by construction $f_2$ is bijective. Its inverse is also given by a Diophantine expression in the identity map. In the same way we can define the function $f_3: \{x \in \Z\mid l\nmid x\} \ra \Z$ by $f_3(x) = f_3(lq+r) = (l-1)q+r-1$. Combining these two maps give the Diophantine and bijective map $f_3^{-1} \circ f_2: \{x \in \Z \mid k\nmid x\} \ra \{x \in \Z \mid l\nmid x\}$. \\
		Using the definition by cases we get that $f:\Z \ra \Z$ defined by
		$$f(x) = \begin{cases} f_1(x) &\text{ if }k\mid x\\ (f_3^{-1}\circ f_2)(x) &\text{ if } k\nmid x\end{cases} $$ 
		is a bijective Diophantine map. In a similar way
		$$f^{-1} = \begin{cases} f_1^{-1}(x) &\text{ if } l\mid x\\ (f_2^{-1}\circ f_3)(x) &\text{ if }l\nmid x\end{cases}$$ 
		is a Diophantine map, so $f$ is a Diophantine automorphism.\\
		We have by construction that $d_l' = f \circ d_k'\circ id$, so $d_k$ and $d_l$ are equivalent.
	\end{example} 
	
	We now extend our definition of similarity and equivalence of models to those that use a Diophantine equivalence map.
	
	\begin{defi} \th\label{diophDefModelsEquiv} 
		Let $d_1:R_1 \ra R_2\modsim_1$ and $d_2:R_1 \ra R_2\modsim_2$ be models of $(R_1,\L_1)$ in $(R_2,\L_2)$. We say that $d_1$ and $d_2$ are \textbf{similar} if there exists a Diophantine automorphism $f$ of $R_2$ such that $\ol{d_1}(R_1) = f(\ol{d_2}(R_1))$ and $\{(x,y)\mid x \sim_1 y\} = \{(f(x),f(y))\mid x\sim_2y\}$. \\
		Suppose that $f_2$ is a Diophantine automorphisms of $R_2$ such that $f_2(x) \sim_1 f_2(y)$ if and only if $x\sim_2y$. Then we define the map $\ol{f_2}: R_2\modsim_2 \ra R_2\modsim_1$ as the unique function such that $\pi_1 \circ f_2 = \ol{f_2} \circ \pi_2$, where $\pi_1: R_2 \ra R_2\modsim_1$ and $\pi_2:R_2 \ra R_2\modsim_2$ are the projections. We say that $d_1$ and $d_2$ are \textbf{equivalent} if there exist Diophantine automorphisms $f_1$ and $f_2$ such that $f_2(x) \sim_1 f_2(y) \LRa x\sim_2 y$ and $d_1 = \ol{f_2} \circ d_2 \circ f_1$. 
	\end{defi} 
	
	\begin{remark}
		We have that Diophantine maps are also Diophantine equivalence maps by using the relation equality as the equivalence relation. If $\sim_1$ and $\sim_2$ are both equality we get back the previous definition, so this is a proper extension. \\
		Also in this setting, $d_1$ and $d_2$ are similar if they are equivalent. We have namely 
		$$\{(x,y)\mid x\sim_1 y\} = \{(f_2(x),f_2(y))\mid f_2(x) \sim_1 f_2(y)\} = \{(f_2(x),f_2(y))\mid x\sim_2 y\}$$ 
		and $f_1$ is bijective, so $\ol{d_1}(R_1) = f_2(\ol{d_1}(f_1(R_1))) = f_2(\ol{d_1}(R_1))$.
	\end{remark} 
	
	\begin{example}
		We will define two models of $\F_2$ in $\Z$. \\
		Define $\sim_1$ by $x\sim_1 y$ if and only if $x\geq 0$ and $y \geq 0$ or $x<0$ and $y<0$. Then define $d_1: \F_1 \ra \Z\modsim_1$ by $d_1(0) = [0]_1 = \N$ and $d_1(1) = [-1]_1 = \Z_{<0}$. Since $\ol{d_1}(1) = \Z_{<0}$, $\ol{d_1}(0) = \N$, $\ol{d_1}(\emptyset) = \emptyset$ and $\ol{d_1}(\F_2) = \Z$ are all Diophantine sets, we get that $d_1$ is a Diophantine equivalence map and therefore a Diophantine model.\\
		Define $\sim_2$ by $x \sim_2 y$ if and only if $2\mid x-y$ and define $d_2: \F_2 \ra \Z\modsim_2$ by $d_2(x) = [x]_2$. The equivalence classes under $\sim_2$ are $\{x \in \Z \mid 2\mid x\}$ and $\{x \in \Z \mid 2\nmid x\}$, which are both Diophantine sets. Since $\emptyset$ and $\Z$ are also Diophantine sets, we get that $d_2$ is a Diophantine equivalence map and thus a Diophantine model. \\
		We define the map $f$ by 
		$$ f(x) = \begin{cases} x & \text{ if } 2\mid x,\ x \geq 0\\
			-x-1 &\text{ if } 2\mid x,\ x < 0\\
			-x-1 &\text{ if } 2\nmid x,\ x\geq 0\\
			x &\text{ if } 2\nmid x,\ x <0.\\ \end{cases} $$
		The function $x \ra -x-1$ is given by a Diophantine expression in the identity function, so it is a Diophantine map by \th\ref{diophd1dioph+d1<Dd2=>d2dioph}. Since $\Z$ is a domain, $\L_R$ is one-sided over $\Z$ and the cases are Diophantine sets, we get with \th\ref{diophmapcases} that $f$ is a Diophantine map. Since $f$ is its own inverse, we get that $f$ is a Diophantine automorphism.\\
		We have that $f([0]_2) = \N = [0]_1$ and $f([1]_2) = \Z_{<0} = [-1]_1$, so it holds that  $f(\ol{d_2}(R_1)) = f(R_2) = \ol{d_1}(R_1)$ and $\{(x,y)\mid x\sim_1 y\} = \{(f(x),f(y))\mid x\sim_2 y\}$. This implies that $f$ is a Diophantine automorphism that shows that $d_1$ and $d_2$ are similar.\\
		We have that $f(x) \sim_1 f(y)$ if and only if $x \sim_2 y$ and $\ol{f} \circ d_2 \circ id = d_1$, so $d_1$ and $d_2$ are also equivalent.  
	\end{example}

	\section{Comparing recursive and Diophantine} 
	
	We start this section by looking at how recursive sets and recursively enumerable sets behave under certain maps. We then compare this behaviour to that of Diophantine sets. The necessary background concerning the theory of recursive functions, recursive sets, recursively enumerable sets and recursive rings can be found in Appendix A of \cite{Shlapentokh}. 
	
	\subsection{Recursive and Recursively Enumerable Sets under Maps}
	
	\begin{assumption}
		All rings will be recursive rings. If they are not finitely generated as a ring or as a field, we assume they come with a fixed recursive representation and all recursive functions and sets are recursive under that specific recursive representation.
	\end{assumption}
	
	\begin{lemma} \th\label{recfrecSRE=>f(S)RE} 
		Let $f: R_1 \ra R_2$ be a recursive function and $S \subset R_1$ a recursively enumerable set. Then $f(S) \subset R_2$ is recursively enumerable. \qed 
	\end{lemma}
	
	\begin{defi}[\th\ref{introDefREpres}]
		We call a function $f: R_1 \ra R_2$ \textbf{recursively enumerable preserving}, or \textbf{r.e.\ preserving} for short, if for all recursively enumerable sets $S \subset R_1$, the set $f(S) \subset R_2$ is recursively enumerable. 
	\end{defi}
	
	\begin{remark}
		Note the similarity between the definition of r.e.\ preserving functions and the definition of a Diophantine map. 
	\end{remark}
	
	\begin{remark} \th\label{recrecsubsetREpres}
		By \th\ref{recfrecSRE=>f(S)RE} we get that a recursive function is always a recursive enumerable preserving function. The converse does not need to hold, see \th\ref{recREpres/=rec}.
	\end{remark}
	
	\begin{lemma} \th\label{recfREpresSrec=>f(S)rec} 
		Let $R_1$ and $R_2$ be rings and $f: R_1 \ra R_2$ a recursive enumerable preserving function. Let $S \subset R_1$ be a recursive set. If $f$ is injective and $f(R_1)$ is a recursive set, then $f(S)$ is a recursive set.
	\end{lemma}
	
	\begin{proof}
		Since $S$ is recursive, the sets $S$ and $S^c = R_1\bs S$ are recursively enumerable. This gives that $f(S)$ and $f(S^c)$ are recursively enumerable. Since $f(R_1)$ is a recursive set, also $R_2\bs f(R_1) = f(R_1)^c$ is a recursively enumerable set. This implies that the set $f(S^c)\cup f(R_1)^c$ is a recursive set. Since $f$ is injective, we get $f(S)^c = f(S^c) \cup f(R_1)^c$, so $f(S)$ is a recursive set. 
	\end{proof}
	
	\begin{remark}
		Both conditions in \th\ref{recfREpresSrec=>f(S)rec} are needed.\\ 
		Suppose that $R_1$ and $R_2$ are infinite. Let $S \subset R_2$ be a recursively enumerable, but not recursive set. There exists a recursive $f: \N \ra R_1$ that enumerates $S$ injectively and there exists a bijective recursive representation $J_1: R_1 \ra \N$. (For an introduction to recursive representations, see \cite{Shlapentokh}.) Then the function $f\circ J_1: R_1 \ra R_2$ is a recursive injective function. The set $R_1$ is recursive, but $(f\circ J_1)(R_1) = f(\N) = S$ is not a recursive set. So only injectivity is not enough.\\
		Define $g: \N \ra R_2$ by 
		$$g(n) = \begin{cases} f(n\div 2) &\text{if }2\mid n\\ n\div2 &\text{if }2\nmid n\\ \end{cases}.$$
		Then $g$ is a recursive function, so the function $g\circ J_1:R_1 \ra R_2$ is also recursive. We have that $g(J_1(R_1)) = g(\N) = \N$, which is a recursive set. On the other hand the set $g(J_1(J_1^{-1}(\{2n\mid n \in \N\}))) = g(\{2n\mid n \in \N\}) = f(\N) = S$ is not recursive, while $J_1^{-1}(\{2n\mid n \in \N\})$ is a recursive set, since the set $\{2n\mid n \in \N\} = \{n\in \N\mid 2\mid n\}$ is recursive. So only asking the image of the ring to be recursive is also not enough.
	\end{remark}
	
	\begin{lemma} \th\label{recfrecSrec=>f-1(S)rec} 
		Let $f: R_1 \ra R_2$ be a recursive function and let $S \subset R_2$ be a recursive set. Then $f^{-1}(S)$ is a recursive set. \qed 
	\end{lemma}
	
	\begin{lemma} \th\label{recfrecinjSRE=>f-1(S)RE} 
		Let $f:R_1 \ra R_2$ be an injective recursive function. Let $S \subset R_2$ be a recursively enumerable set. Then $f^{-1}(S) \subset R_1$ is a recursively enumerable set. \qed 
	\end{lemma}

	\begin{lemma} \th\label{recREpres/=rec} 
		Let $R_1,R_2$ be infinite rings. Then there exists a surjective r.e.\ preserving function $f: R_1 \ra R_2$ that is not a recursive function.
	\end{lemma}
	
	\begin{proof}
		We start by constructing a surjective r.e.\ preserving function $g: \N \ra \N$ that is not a recursive function. There exists a recursively enumerable set $S \subset \N$ that is not recursive. Let $h: \N \ra \N$ recursively enumerate this set. The set $S$ is infinite since otherwise it would be recursive, so we can choose $h$ to be injective. Define $g: \N \ra \N$ by 
		$$g(n) = \begin{cases} 0 &\If n\notin S\\ 1+(\mu m)(h(m)=n) &\If n \in S\end{cases}.$$
		Let $A \subset \N$ be a recursively enumerable set. If $A \cap S^c = \emptyset$, we get $g(A) = \{1+(\mu m)(h(m)=n)\mid n \in A\}$. If $A$ is empty, so is $g(A)$, so then it is recursively enumerable. If not, let $f_A$ recursively enumerate $A$. Then $g_A(n) =1+(\mu m)(h(m)=f_A(n))$ is a recursive function that enumerates $g(A)$. If $A \cap S^c \not= \emptyset$, then $g(A) = g(A \cap S) \cup \{0\}$. Since $S$ is recursively enumerable, so is $A \cap S$. So by the previous part we get that $g(A)$ is a recursively enumerable set. This implies that $g$ is r.e.\ preserving. \\
		We have that $\{0\}$ is a recursive set. On the other hand, $g^{-1}(0) = S^c$ is not a recursively enumerable set, so not a recursive set. By \th\ref{recfrecSrec=>f-1(S)rec} this gives that $g$ is not a recursive function. \\
		Since $h$ is an injective function, we can also write $g$ as 
		$$g(n) = \begin{cases} 0 &\If n \notin S\\ 1+ h^{-1}(n)&\If n \in S\end{cases}.$$
		Since $h(\N) = S$, this gives that $g$ is a surjective function.\\
		Now let $J_1$ and $J_2$ be bijective recursive representations of $R_1$ respectively $R_2$. Define $f = J_2^{-1} \circ g \circ J_1$. Let $A \subset R_1$ be a recursively enumerable set. Then by \th\ref{recfrecSRE=>f(S)RE} we get that $J_1(A)$, hence $g(J_1(A))$, hence also $(J_2^{-1} \circ g \circ J_1)(A)$ is recursively enumerable. So $J_2^{-1} \circ g \circ J_1$ is an r.e.\ preserving function. We have that $f$ is recursive if and only if its translation is the restriction of a recursive function. Now $J_2 \circ f \circ J_1^{-1} = J_2 \circ J_2^{-1} \circ g \circ J_1 \circ J_1^{-1} = g|_{J_1(R_1)} = g$ is not the restriction of a recursive function, since $J_1(R_1) = \N$. So $f: R_1 \ra R_2$ is an r.e.\ preserving, but not a recursive function. Since $J_1^{-1}$, $g$ and $J_2$ are surjective, so is $f$. 
	\end{proof}

	\subsection{Recursive versus Diophantine} \label{parRecvsDioph} \phantom{=}\\
	In this paragraph, we will look at the connections between recursive functions, recursively enumerable functions, Diophantine functions and functions that are given by a Diophantine expression in some kind of identity function.
	
	\begin{defi}
		Let $(R,\L)$ be a ring with a language. We call the interpretation of $\L$ in $R$ \textbf{recursive} if for all function $f \in \L$, the interpretation $f^R$ is a recursive function and for all relations $S \in \L$, the interpretation $S^R$ is a recursive set.   
	\end{defi}
	
	\begin{lemma} \th\label{recdiophset=>REset}
		Let $(R,\L)$ be a recursive ring with a recursive interpretation of $\L$. Let $S \subset R^k$ be a Diophantine set. Then $S$ is a recursively enumerable set. \qed 
	\end{lemma}
	
	\begin{lemma} \th\label{recRLrecRkLkrec} 
		Let $(R,\L)$ be a ring with a recursive interpretation of $\L$, then $(R^k,\L^{(k)})$ is a ring with a recursive interpretation of $\L^{(k)}$. \qed 
	\end{lemma}
	
	\begin{assumption}
		For the rest of this chapter, $(R,\L)$ will be a recursive ring $R$ with a recursive interpretation of the language $\L$.
	\end{assumption}
	
	The solution of Hilbert's 10th problem over $\Z$ or over $\N$ uses the equivalence between recursively enumerable and Diophantine sets. The theorem we need is the following variant.
	
	\begin{thm} \th\label{recNrecfunc=diophexprfunc}
		Let $f: \N^k \ra \N$ be a recursive function. Then there exists $n \in \N$ and polynomials $P_1,\ldots,P_n \in \Z[X_1,\ldots,X_{l+k+1}]$ such that
		$$f(\vec{x})=z \LRa \ex \vec{y} \in \Z^l,\ \va 1\leq i \leq n,\ P_i(\vec{x},\vec{y},z)=0.$$
		This is equivalent to $f$ being given by a Diophantine expression in $id: \Z \ra \Z$ where we give $\Z$ the language of rings. 
	\end{thm}
	
	\begin{proof}
		This is a reformulation of Theorem 5.24 of \cite{H10PMurtyFodden}.
	\end{proof}
	
	\begin{cor} \th\label{recNkrecfunc=diophexprfunc}  
		Let $f: \N^k \ra \N^l$ be a recursive function. Then there exists $n \in \N$ and polynomials $P_1,\ldots,P_n \in \Z[X_1,\ldots,X_{k+l+m}]$ such that 
		$$f(\vec{x})  = \vec{z} \LRa \ex \vec{y} \in \Z^m,\ \va 1\leq i \leq n,\ P_i(\vec{x},\vec{y},\vec{z}) = 0.$$
		Again this is equivalent to $f$ being given by a Diophantine expression in $id: \Z \ra \Z$ with the language of rings on $\Z$.
	\end{cor}
	
	\begin{thm}[MRDP] \th\label{recMRDPNk}
		Let $S \subset \N^k$ be a recursively enumerable set. If we give $\N$ the language of rings, then $S$ is a Diophantine set.
	\end{thm}
	
	\begin{remark}
		The version with $S \subset \N$ is the theorem that Yuri Matiyasevich, Julia Robinson, Martin Davis and Hilary Putnam proved. The theorem is named MRDP after them. Note that sometimes MRDP is used for saying the Hilbert's tenth problem has a negative answer, which is weaker then recursively enumerable and Diophantine sets being the same. 
	\end{remark}
	
	\begin{defi}
		We say that \textbf{MRDP holds} in $(R,\L)$ if for all $S \subset R$ recursively enumerable, $S$ is a Diophantine set in $(R,\L)$. 
	\end{defi}

	\begin{lemma} \th\label{recMRDPRfinite}
		Let $(R,\L)$ be a finite integral domain with a language $\L$ that consists of the language of rings with possible extra constants. If $R = S_c'(\L)$, then MRDP holds in $(R,\L)$. \qed 
	\end{lemma}
	
	\begin{lemma} \th\label{recdiophmap=>REpresmapifMRDPR1}
		Let $(R_1,\L_1)$ and $(R_2,\L_2)$ be rings with a language. Suppose that MRDP holds in $(R_1,\L_1)$. Let $d: R_1 \ra R_2$ be a Diophantine map, then $d$ is r.e.\ preserving. \qed 
	\end{lemma}
	
	\begin{lemma} \th\label{recREpresmap=>diophmapifMRDPR2} 
		Let $(R_1,\L_1)$ and $(R_2,\L_2)$ be rings with a language. If MRDP holds in $(R_2,\L_2)$, an r.e.\ preserving map $f:R_1 \ra R_2$ is also a Diophantine map. \qed 
	\end{lemma}
	
	\begin{prop} \th\label{recMRDPR=>MRDPRk} 
		Let $(R,\L)$ be an infinite ring with a language. If MRDP holds in $(R,\L)$, then it holds in $(R^k,\L^{(k)})$.
	\end{prop}
	
	\begin{proof}
		Let $J: R \ra \N$ be a bijective recursive representation of $R$. We have that $J^{-1}$ is a recursive function. Since MRDP holds in $(R,\L)$, we get with \th\ref{recREpresmap=>diophmapifMRDPR2} that $J^{-1}$ is a Diophantine map. Since $(J^{-1})^k$ is given by a Diophantine expression in the injective Dipohantine map $J^{-1}$, we get with \th\ref{diophd1dioph+d1<Dd2=>d2dioph} that $(J^{-1})^k$ is a Diophantine map. 
		Now let $S \subset R^k$ be recursively enumerable. Then by definition, $J^k(S)$ is a recursively enumerable set as $J^k:R^k \ra \N^k$ is a recursive representation of $R^k$. Since MRDP holds in $\N^k$, we get that $J^k(S)$ is a Diophantine set. Then $S = (J^{-1})^k (J^k(S))$ is a Diophantine set, as $(J^{-1})^k$ is a Diophantine map. So MRDP holds in $(R^k,\L^{(k)})$.  
	\end{proof}
	
	\begin{lemma} \th\label{recMRDPRk=>MRDPR} 
		Let $(R,\L)$ be a ring with a language. If there exists a $k \in \Z_{>0}$ such that MRDP holds in $(R^k,\L^{(k)})$, then MRDP holds in $(R,\L)$. \qed 
	\end{lemma}
	
	\begin{lemma} \th\label{recrecbij} 
		Let $R_1$ and $R_2$ be infinite recursive rings. Then there exists a recursive bijection $\phi: R_1 \ra R_2$. 
	\end{lemma}
	
	\begin{proof} 
		Let $J_1$ and $J_2$ be bijective recursive representations of $R_1$ respectively $R_2$. Then the function $\phi: R_1 \ra R_2$ defined by $\phi = J_2^{-1} \circ J_1$ is a recursive bijection.
	\end{proof} 
	
	\begin{prop}
		Let $(R_1,\L_1)$ and $(R_2,\L_2)$ be rings and let $\phi:R_1 \ra R_2$ be an injective recursive function. Suppose that MRDP holds in $R_2$. Then MRDP holds in $R_1$ if and only if $\phi^{-1}:\phi(R_1) \ra R_1$ is a Diophantine map. \qed 
	\end{prop}
	
	\begin{remark}
		Note that $\phi(R_1)$ is recursive, so $\phi(R_1)$ and $\phi(R_1)^{c}$ are recursively enumerable. Since MRDP holds in $R_2$, this means that they are both Diophantine sets, so we can speak of a Diophantine map from $\phi(R_1)$. See also \th\ref{diophmapfromsubset}.
	\end{remark}
	
	\begin{lemma} \th\label{recdrec+MRDPR2=>ddiophexprphi} 
		Let $(R_1,\L_1)$ and $(R_2,\L_2)$ be rings with a language. Let $d: R_1 \ra R_2$ be a recursive function and let $\phi:R_1 \ra R_2$ be a recursive and bijective function. If MRDP holds in $(R_2,\L_2)$, then $d$ is given by a Diophantine expression in $\phi$.
	\end{lemma}
	
	\begin{proof}
		The map $(\phi,d): R_1 \ra R_2^2$ defined by $(\phi,d)(x) = (\phi(x),d(x))$ is a recursive function. Since $R_1$ is a recursively enumerable set, this gives that $S = \{(\phi(x),d(x))\mid x \in R_1\} \subset R_2^2$ is a recursively enumerable set (\th\ref{recfrecSRE=>f(S)RE}). \\
		By \th\ref{recMRDPR=>MRDPRk} we get that MRDP also holds in $(R_2^2,\L_2^{(2)})$. This gives that $S$ is a Diophantine set, so there exist basic Diophantine $\L_2$-formulas $f_1,\ldots,f_n$ such that
		$$S = \{(a,b) \in R_2^2\mid \ex \vec{v} \in R_2^k,\ \va 1\leq i \leq n,\ f_i(a,b,\vec{v})\}.$$
		Hence
		$$d(x)=b \LRa \ex \vec{v} \in R_2^k,\ \va 1\leq i \leq n,\ f_i(\phi(x),b,\vec{v}),$$
		so $d$ is given by a Diophantine expression in $\phi$.
	\end{proof}
	
	\begin{lemma} \th\label{recddiophphi=>drec} 
		Let $(R_1,\L_1)$ and $(R_2,\L_2)$ be rings with a language. Let $\phi: R_1 \ra R_2$ be a bijective recursive function. Suppose that $d: R_1 \ra R_2$ is given by a Diophantine expression in $\phi$. Then $d$ is a recursive function. \qed 
	\end{lemma}
	
	\begin{remark}
		We can summarize the above results in a diagram, see the end of the introduction.
	\end{remark} 

	\begin{lemma} \th\label{recrecinjR1finiteR2} 
		Let $R_1$ be a finite ring and $R_2$ an infinite ring. Then there exists an injective recursive function $\phi: R_1 \ra R_2$.
	\end{lemma}
	
	\begin{lemma}
		Let $f:R_1 \ra R_2$ be any function and suppose that $R_1$ is finite. Then $f$ is a recursive function. \qed 
	\end{lemma}
	
	\begin{lemma}
		Let $R_1$ be a finite ring with language $\L_1$ and let $R_2$ be an integral domain with language $\L_2$. Suppose that $\L_2$ is the language of rings with set of constants $S_c$ such that $S_c'(\L_2) = R_2$, so all singletons are Diophantine. Let $\phi: R_1 \ra R_2$ be an injective map. Then $f:R_1 \ra R_2$ is given by a Diophantine expression in $\phi$. Further, $f$ and $\phi$ are Diophantine maps. \qed 
	\end{lemma}
	
	\begin{remark}
		Let $(R_1,\L_1)$ and $(R_2,\L_2)$ be rings such that $R_1$ is finite and $|R_2| \geq |R_1|$. Suppose further that $\L_2$ is the language of rings with possible extra constants and $R_2$ is an integral domain. Let $\phi:R_1 \ra R_2$ be an injective map. Then we get that:
		\begin{align*} 
			\{d:R_1 \ra R_2\} &= \{d: R_1 \ra R_2\mid d \text{ Diophantine}\}\\ 
			&= \{d:R_1 \ra R_2\mid d\text{ r.e.\ preserving}\}\\
			&= \{d:R_1 \ra R_2\mid d\text{ recursive}\}\\
			&= \{d:R_1 \ra R_2\mid d\text{ given by a Diophantine expression in }\phi\}
		\end{align*}
		If $|R_1|>|R_2|$, we have no injective $\phi:R_1 \ra R_2$ to express our maps in, but the other equalities still hold. 
	\end{remark}
	
	\begin{lemma} \th\label{recR2finitefrecpres} 
		Let $f:R_1 \ra R_2$ be a map and suppose that $R_2$ is finite. Then $f$ is r.e.\ preserving. \qed 
	\end{lemma}
	
	\begin{lemma}
		Let $R_1$ be an infinite ring and let $R_2$ be a finite non-zero ring. Then there exists a function $f:R_1 \ra R_2$ that is r.e.\ preserving, but not recursive. \qed 
	\end{lemma}
	
	\begin{remark}
		Let $R_1$ be an infinite ring and $R_2$ be a finite non-zero ring. Since $|R_2|\leq |R_1|$, there is no injective function from $R_1\ra R_2$. This means we have no $\phi$ to make sense of the set $\{f:R_1 \ra R_2\mid \phi \leq_D f\}$, so we omit this from the diagram.\\
		Let $\L_2$ be the language on $R_2$. If $\L_2$ is the language of rings with possible extra constants and $R_2$ is an integral domain, we can use \th\ref{diophmapfinite}. Then we can summarize the situation:
		$$\{f:R_1 \ra R_2\mid f\text{ recursive}\} \subsetneq \{f: R_1 \ra R_2 \mid f\text{ r.e.\ preserving}\} = \{d:R_1 \ra R_2\mid d\text{ Diophantine}\}.$$ 
	\end{remark}

	\section{Hilbert's Tenth Problem} 
	
	In this section, we first formalize Hilbert's tenth problem. Then we take a look into how Diophantine maps can be used to transfer a negative or positive answer to Hilbert's tenth problem from one ring to another.
	
	\subsection{A Formal Version of Hilbert's Tenth Problem} 
	
	\begin{notation}
		Let $(R,\L)$ be a set with a language and let $t$ be an $\L$-term. Let $x$ be a variable of $\L$ and let $y$ be a variable or an element of $S_c'(\L) = \{r \in R \mid \{r\}\text{ is a Diophantine set}\}$. We denote by $t[x/y]$ the term of $\L$ where all occurrences of $x$ are replaced by a $y$.
	\end{notation}
	
	\begin{lemma} \th\label{H10PnumberingLterms} 
		Let $(R,\L)$ be a set with a countable language $\L$. Let $R_0 \subset S_c'(\L)$ be a ring. Suppose that we have a recursive representation $j: R_0 \ra \N$. Let $\numL$ be a function that numbers the functions and constants of $\L$. Furthermore, define the function $\sym:\{\text{constants and functions of }\L\}\sqcup \N \sqcup R_0 \ra \N$ by $\sym(l) = 2^{\numL(l)}$ for $l \in \L$, $\sym(n) = 3^n$ for $n \in \N$ and $\sym(r) = 5^{j(r)}$ for $r \in R_0$. Then there exists a numbering $J: \{\L\text{-terms}\} \ra \N$ and a recursive function $\ev: R_0 \ti \N^2 \ra \N$ such that $\ev(a,n,J(t)) = J(t[x_n/a])$, where $\{x_i\mid i \in \N\}$ are the variables used. 
	\end{lemma}
	
	\begin{proof}
		We start by writing the terms of $\L$ as a string of symbols in the formal way, so we do not need brackets, see also Paragraph 2.2 of \cite{SetsModelsProofs}. Then define 
		$$J(s_1,\ldots,s_n) = \prod_{i=1}^n p_i^{\sym(s_i)}.$$
		It is straightforward to show that $J$ together with the function $\ev$ that replaces all $p_i^{3^n}$ by $p_i^{5^{j(a)}}$ satisfies the indicated properties.
	\end{proof}
	
	\begin{lemma} \th\label{H10Pev} 
		Let $(R,\L)$, $R_0$ and $J$ be as in \th\ref{H10PnumberingLterms}. Then there exist recursive functions $\ev_k: R_0^k \ti \N \ra \N$ such that $\ev_k(\vec{a},J(t)) = J(t[x_1/a_1][x_2/a_2]\ldots [x_k/a_k])$. \qed 
	\end{lemma} 
	
	Since algorithms can not work with an uncountable amount of possibilities, the set $R$ and language $\L$ need to satisfy some conditions for Hilbert's tenth problem to make sense. This leads to the following definition.
	
	\begin{defi} \th\label{H10Pdefcompatibility}
		Let $(R,\L)$ be a set with a language. Suppose that $R_0$ is a recursive ring such that $S_c(\L) \subset R_0 \subset S_c'(\L)$. Define for all functions $f \in \L$ the relation $S_f = \{(x_1,\ldots,x_r,f(\vec{x})) \mid \vec{x} \in R^r\}$. Let $\L'$ be the language on $R_0$ that has for each $c \in \L$ a constant $c$, for each function $f \in \L$ the relation $S_f \cap R_0^r$ and for each relation $S$ in $\L$ the relation $S \cap R_0^r$.\\
		We say that $R_0$ is \textbf{compatible} with $(R,\L)$ if $S_c'(\L') = R_0$ and if the interpretation of $\L'$ in $R_0$ is recursive. 
	\end{defi}
	
	\begin{defi} \th\label{H10PexistenceH10P}
		Let $(R,\L)$ be a set with a language. We say that \textbf{Hilbert's tenth problem exists over $(R,\L)$} if $\L$ is countable and there exists a recursive ring $R_0$ such that $S_c(\L) \subset R_0 \subset S_c'(\L)$ and $R_0$ is compatible with $(R,\L)$. 
	\end{defi}
	
	\begin{lemma} 
		If $(R,\L)$ is a set with a countable language such that $S_c(\L)$ is finite, then Hilbert's tenth problem exists. 
	\end{lemma} 
	
	\begin{proof} 
		We take $R_0 = S_c(\L)$ and equip it with an arbitrary ring structure. Then $R_0$ is a finite ring, so it is an recursive ring. Also all finite sets are recursive, so $\L'$ will have a recursive interpretation in $R_0$. This gives that $R_0$ is compatible with $(R,\L)$, so Hilbert's tenth problem exists over $(R,\L)$, since $\L$ is countable. 
	\end{proof}
	
	\begin{remark}
		The following definition is a generalization of the formalization of Hilbert's tenth problem in \cite{Shlapentokh}.
	\end{remark}
	
	\begin{defi}
		Let $R$ be a set and let $\L$ be a countable language with an interpretation in $R$. Suppose that $j$ is a recursive representation of a compatible ring $R_0$ with $S_c(\L) \subset R_0 \subset S_c'(\L)$. Let $J$ be the numbering of \th\ref{H10PnumberingLterms}. Let $\{S_i \subset R^{k_i} \mid i \in I \subset \N\}$ be the relations of $\L$, including equality. \\
		We say that the \textbf{strict version of Hilbert's tenth problem over $R$ has a positive answer} if for all $i \in I$ there exists a recursive function $f_i: \N^{k_i} \ra \{0,1\}$ such that for all $\L$-terms $t_1,\ldots,t_{k_i}$ we have that $f_i(J(t_1),\ldots,J(t_{k_i})) = 1$ if and only if there exist $x_1,\ldots,x_m \in R$ such that $(t_1(\vec{x}),\ldots,t_{k_i}(\vec{x})) \in S_i$. \\ 
		The \textbf{strict version of Hilbert's tenth problem over $R$} has a negative answer if for all choices of $R_0$ and all recursive representations of $R_0$ such a function $f$ does not exist.\\
		We say that \textbf{Hilbert's tenth problem over $R$ has a positive answer} if there exists a recursive function $f: \bigcup_{n=0}^\infty \N^n \ra \{0,1\}$ such that for all $\L$-terms $t_{11},\ldots,t_{lk_{i_l}}$ we have that
		$$f(l,i_1,J(t_{11}),\ldots,J(t_{1k_{i_1}}),\ldots,i_l,J(t_{l1}),\ldots,J(t_{lk_{i_l}})) = 1$$
		if and only if the system 
		$$(t_{11},\ldots,t_{1k_{i_1}})\in S_{i_1},\ldots,(t_{l1},\ldots,t_{lk_{i_l}}) \in S_{i_l}$$ 
		has a solution over $R$. \\
		\textbf{Hilbert's tenth problem over $R$ has a negative answer} if for all choices of $R_0$ and for all recursive representations of $R_0$ such an $f$ does not exists. 
	\end{defi}
	
	\begin{remark}
		The choice of $R_0$ does not matter for the general version of Hilbert's tenth problem, since we can replace a coefficient in $S_c'(\L)$ by a set of Diophantine equations with coefficients in $S_c(\L) \subset R_0$.
	\end{remark}
	
	\begin{thm} \th\label{H10PchoicerecreprR0}
		Let $(R,\L)$ be a set with a finite language and let $S_c(\L) \subset R_0 \subset S_c'(\L)$ be a set. Suppose that $(R_0,+_1,\cdot_1)$ with $j_1$ and $(R_0,+_2,\cdot_2)$ with $j_2$ are two ring structures with a recursive representation that are compatible with $(R,\L)$. If Hilbert's tenth problem has a positive answer with respect to $(R_0,j_1)$, then Hilbert's tenth problem has a positive answer with respect to $(R_0,j_2)$. In other words, in this case it does not matter which recursive representation on $R_0$ we take. 
	\end{thm}
	
	For the proof of this theorem we need the following two results.
	
	\begin{lemma} \th\label{H10Puitpakkenterms} 
		Let $(R,\L)$ be a set with a language. Let $R_0$ be a recursive ring with the property that $S_c(\L) \subset R_0 \subset S_c'(\L)$. Let $\numL$, $\sym$ and $J$ be the numbering of the constants and functions of $\L$, the symbols in an $\L$-term and the $\L$-terms of \th\ref{H10PnumberingLterms}, respectively.\\
		Suppose that there exist a recursive function $F: \N \ra \N$ such that $F(\numL(c)) = 0$ for all constants of $\L$ and $F(\numL(f)) = r$ for all $r$-ary functions $f$ of $\L$. \\
		There exist a recursive function $u: \N \ra \N \ti \bigcup_{n=0}^\infty \N^n$ such that $u(J(x_n)) = (\sym(x_n))$, $u(J(c)) = (\sym(c))$ and $u(J(f(t_1,\ldots,t_n))) = (\sym(f),J(t_1),\ldots,J(t_n))$. 
	\end{lemma} 
	
	\begin{proof}
		We will describe the algorithm that computes $u$.\\
		Suppose that we have input $t$. First set $i$ equal to $l(t)$, the largest prime factor of $t$. Let $S \subset \bigcup_{i=0}^\infty \N^n$ be the empty sequence. \\
		Write $\rho(i,t)$ for the function that gives the exponent of $p_i$ in the prime factorization of $t$. Then while $i>1$ we do the following:
		\begin{itemize}
			\item If $\rho(i,t) = 3^n$, so the $i$-th symbol of $t$ is the variable $x_n$, add $2^{\rho(i,t)}$ to $S$.
			\item If $\rho(i,t) = 5^n$, so the $i$-th symbol of $t$ is an element of $R_0$, add $2^{\rho(i,t)}$ to $S$. 
			\item If $\rho(i,t) = 2^n$ and $F(n)=0$, we have that the $i$-th symbol of $x$ is a constant. Also in this case we add $2^{\rho(i,t)}$ to $S$.
			\item If $\rho(i,t) = 2^n$ and $F(n) = r>0$, we have that the $i$-th symbol of $x$ is an $r$-ary function $f$. Let $t_1,\ldots,t_r$ be the final $r$ numbers of $S$, with $t_1$ being the last. Remove these numbers from $S$. We want to add the term $J(f(J^{-1}(t_1),\ldots,J^{-1}(t_r)))$ to $S$. \\
			We can compute this in the following way. Initialize the counter $m=1$ and the result $N=1$. \\
			For each $1 \leq k \leq r$ and for each $1\leq j \leq l(k)$, multiply the result $N$ by $p_m^{\rho(j,t_k)}$ and add one to the counter $m$. When done we have that $N = J(f(J^{-1}(t_1),\ldots,J^{-1}(t_r)))$, so we add $N$ to $S$.
		\end{itemize}
		When $i=1$, we can output $(\rho(1,x),S)$, since $S$ will be empty if the first symbol of $x$ is not a function and $S$ will be of the form $(J(t_1),\ldots,J(t_r))$ if the first symbol of $x$ is an $r$-ary function. \\
		Since all functions used are recursive functions, $u$ is also a recursive function. 
	\end{proof} 
	
	\begin{prop} \th\label{H10Puniquerecrepr} 
		Let $(R_1,j_1)$ and $(R_2,j_2)$ be rings with a recursive representation such that $R_1 = R_2$ as sets. Let $\L$ be a finite language that has a recursive interpretation in $R_1$ and $R_2$. Suppose that $\{x\} \subset R_i$ is a Diophantine set for all $x \in R_i$. Then there exists a recursive function $\phi: \N \ra \N$ such that $\phi \circ j_1 = j_2$. 
	\end{prop}
	
	\begin{proof}
		Let $x\in R_1$. Since $\{x\}$ is a Diophantine set, there exists an expression $$F_x(x',\vec{y}) = (g_{11}(x',\vec{y}),\ldots,g_{1i_1}(x',\vec{y})) \in S_{j_1} \wedge \ldots \wedge (g_{n1}(x',\vec{y}),\ldots,g_{ni_n}(x',\vec{y})) \in S_{j_n}$$
		with $g_{ij}$ terms of $\L$ and $S_i$ relations of $\L$ such that $x$ is the only $x'$ for which there exists a $\vec{y}$ such that $F_x(x',\vec{y})$ is true. Furthermore, we can require that the $g_{ij}$ do not contain elements of $R_1 = R_2$, since we can rewrite those with extra equations. This gives that $x = \pi_1(\mu(x',\vec{y})(F_x(x',\vec{y})))$.\\
		For $m \in j_1(R_1)$ we define 
		$$\phi(m) = \pi_1( \mu(x',\vec{y})(j_2(F_{j_1^{-1}(m)})(j_2(y_1),\ldots,j_2(y_k))))$$
		where $j_2(F_x)$ is the translation of $F_x$ that uses $j_2(c)$ for a constant $c$, $j_2(S_i)$ for a relation $S_i$ and the translation $T_{2,f}(a_1,\ldots,a_r) = j_2(f(j_2^{-1}(a_1),\ldots,j_2^{-1}(a_r)))$ for a function $f$. If $m \notin J_1(R_1)$ we set $\phi(m) = 0$. \\
		We get that
		\begin{align*} 
			\phi(j_1(x)) &= \pi_1(\mu(x',\vec{y})(j_2(F_x)(j_2(x'),\ldots,j_2(y_1),\ldots,j_2(y_m)))) \\
			&= \pi_1(\mu(x',\vec{y})(j_2(F_x(x',\vec{y}))))  = j_2(x). 
		\end{align*} 
		This proves that $\phi$ is the function we need and we are left with proving that it is a recursive function. \\
		We have the functions $J_1$ and $J_2$ from \th\ref{H10PnumberingLterms} that number the $\L$-terms, one with respect to $j_1$ and one with respect to $j_2$. This gives that the set of $\L$-formulas can be seen as a subset of $\cup_{n=1}^\infty \N^n$ in the same way as in the definition of Hilbert's tenth problem. This implies that it makes sense to speak of recursive functions from the set of $\L$-formulas and it also implies that we can minimize over $\L$-formulas. \\
		Note that we have 2 embeddings, one using $j_1$ and one using $j_2$, but they are the same on the formulas that do not use elements of $R_0$, but only constants, variables and functions, since then $J_1$ and $J_2$ agree. \\
		Since $j_1^{-1}(m)$ is unique, we can rewrite $\phi$ to
		$$\phi(m) = \pi_1( \mu(x',\vec{y'})(\mu (F)(j_1(\pi_1(\mu(x,\vec{y})(F(x,\vec{y}))) = m))(j_2(y_1),\ldots,j_2(y_k)))),$$
		where $F$ runs over all $\L$-formulas that do not use elements of $R_1$ or $R_2$.\\
		In the following we will define functions from parts of the $\L$-formulas to $R_i$. We will use the embedding that uses $j_i$ if we map to $R_i$. \\
		We can define the function $\EvTerms_1: \{\L$-terms without variables $\} \ra R_1$ as the function that takes a term without variables and output the element of $R_1$ it evaluates to. We can define it inductively using the function $u$ from \th\ref{H10Puitpakkenterms}. Since $\L$ is finite, we have a finite number of base cases and a finite number of case distinctions in the induction step. Since $\L$ has a recursive interpretation in $R_1$, all steps involved use recursive functions. This implies that $\EvTerms_1$ is a recursive function. In the same way we have a recursive function $\EvTerms_2$ that use $j_2$.\\
		Let $\ev_i:R_i \ti \N^2 \ra \N$ be the function of \th\ref{H10PnumberingLterms} such that $\ev_i(a,n,J_i(t)) = J_i(t[x_n/a])$. \\
		With the functions $\ev_i$ we can make the partial function $\EvTermsVar_i$ from $\{\L\text{-terms}\}\ti \cup_{n=0}^\infty \N^n$ to $R_i$ by setting $\EvTermsVar_i(J_i(t),a_1,\ldots,a_n) = \EvTerms_i(J_i(t[x_1/a_1,\ldots,x_n/a_n]))$ if $J_i(t)$ contains $n$ or less variables and make it undefined otherwise. This is a recursive function since strong induction is recursive and $\EvTermsVar_i(J_i(t),a_1,\ldots,a_n) = \EvTermsVar_i(\ev_i(n,a_n,J_i(t)),a_1,\ldots,a_{n-1})$. \\
		Now we create the final set of partial recursive functions $\Ev_i: \{\L\text{-formulas}\} \ti \cup_{n=0}^\infty R_i^n$ that evaluate the formula $F$ in $(a_1,\ldots,a_n)$ if $F$ uses $n$ or less variables and is otherwise undefined. These functions are recursive since 
		\begin{align*}
			&\phantom{=}\Ev_i((m,k_1,J_i(t_{11}),\ldots,J_i(t_{1l_1}),\ldots,k_n,J_i(t_{m1}),\ldots,J_i(t_{ml_m})),a_1,\ldots,a_n) \\
			&= \prod_{j=1}^m \mathbbm{1}_{k_j}(\EvTermsVar_i(J_i(t_{j1}),a_1,\ldots,a_n),\ldots,\EvTermsVar_i(J_i(t_{jl_j}),a_1,\ldots,a_n)),
		\end{align*} 
		where $\mathbbm{1}_{k_j}$ is the indicator function of the relation $S_{k_j}$. 
		With this we can rewrite $\phi$ to
		$$\phi(m) = \pi_1( \mu(x',\vec{y'})\Ev_2((\mu (F)(j_1(\pi_1(\mu(x,\vec{y})(\Ev_1(F,x,\vec{y}))) = m))),y'_1,\ldots,y'_k)).$$
		Note that $F$ runs over the $\L$-formulas that do not use elements of $R_1$ or $R_2$, so we can ``switch" embeddings of the $\L$-formulas in between using $\Ev_1$ and $\Ev_2$. All functions used here are recursive functions, so $\phi$ is recursive.
	\end{proof}
	
	\begin{proof}[Proof of \th\ref{H10PchoicerecreprR0}]
		With \th\ref{H10Puniquerecrepr} there exists a recursive function $\phi: \N \ra \N$ such that $\phi \circ j_1= j_2$, since $\L$ is finite and compatible with $R_1$ and $R_2$. \\
		Let $J_i: \{\L\text{-terms}\} \ra \N$ be the functions of \th\ref{H10PnumberingLterms} for the recursive representations $j_1$ and $j_2$. Then there exists a recursive function $\psi: \N \ra \N$ such that $\psi \circ J_1 = J_2$, since we need to replace all $p_i^{5^{j_1(a)}}$ by $p_i^{5^{j_2(a)}} = p_i^{5^{\phi(j_1(a))}}$, which can be done recursively. \\
		Then we can define the recursive function $T: \{\L\text{-formulas}\} \ra \{\L\text{-formulas}\}$ by
		\begin{align*} 
			&T(n,i_1,J_1(t_{11}),\ldots,J_1(t_{1m_1}),\ldots,i_n,J_1(t_{n1}),\ldots,J_1(t_{nm_n}))\\
			&=(n,i_1,\psi(J_1(t_{11})),\ldots,\psi(J_1(t_{1m_1})),\ldots,i_n,\psi(J_1(t_{n1})),\ldots,\psi(J_1(t_{nm_n})))
		\end{align*} 
		Note that $J_2$ is injective, so if we restrict $\psi$ to the domain of $J_1$, it is injective. This means we can invert it and extend it to a recursive function $\psi': \N \ra \N$. Then we have $J_1 = \psi'\circ \psi \circ J_1 = \psi '\circ J_2$. This implies that $T$ has an inverse by symmetry and is therefore a bijection.\\
		If now $f:\{\L\text{-formulas}\} \ra \{0,1\}$ is the function that decides Hilbert's tenth problem over $(R,\L)$ with respect to $(R_0,j_1)$, then $f \circ T^{-1}$ decides Hilbert's tenth problem with respect to $(R_0,j_2)$. 
	\end{proof} 
	
	\begin{remark}
		Suppose that $R$ is an integral domain, $\L_R \subset \L$ and $\L$ contains no other relations than equality. If further $\Frac(S_c'(\L))$ is not algebraically closed, then by \th\ref{diophsystem=singlecor} we get that solving a single Diophantine equation is equivalent to solving a system of Diophantine equations. This means that the strict version of Hilbert's tenth problem and the problem itself are equivalent in this case. Note that $(\Z,\L_R)$ fulfills these conditions.  
	\end{remark} 
	
	\begin{remark}
		Note that if Hilbert's tenth problem over $(R,\L)$ has a positive answer, then so has the strict version. If however Hilbert's tenth problem has a negative answer, we do not have a proof of or a counterexample to that the strict version also has a negative answer. It does hold that if the strict version has a negative answer over $(R,\L)$, then the normal version has as well.
	\end{remark} 
	
	\begin{remark}
		The choice of language is very important. We can take the ring $\F_q[T]$ with the language $\{0,1,t,+,\cdot\}$. Then Hilbert's tenth problem has a negative answer by work of Denef, \cite{DenefRTcharp}. Suppose we leave out $T$ and take the language of rings on $\F_q[T]$. Then a (system of) Diophantine equation(s) in this language are polynomials with coefficients in $\F_q$. This means it has a solution over $\F_q[T]$ if and only if it has a solution over $\F_q$, since we can substitute $0$ for $T$. So Hilbert's tenth problem over $(\F_q[T],\L_R)$ has a positive answer.  
	\end{remark}
	
	\begin{lemma}
		Let $(R,\L)$ be a set with a countable language such that $S_c'(\L) = R$ and $R$ has a recursive representation $j$. If Hilbert's tenth problem over $R$ has a positive answer (under $j$), then all Diophantine subsets of $R^k$ are recursive sets (under $j$).
	\end{lemma} 
	
	\begin{proof}
		Let $S \subset R^k$ be a Diophantine set. Then we can write
		$$S = \{\vec{x} \in R^k \mid \ex \vec{y} \in R^{k'},\ \va 1\leq j \leq m,\ (t_{j1},\ldots,t_{jm_j}) \in S_j \}$$
		where the $t_{ij}$ are $\L$-terms and the $S_j$ are (not necessary distinct) relations of $\L$, which can include equality. Let $J$ be as in \th\ref{H10PnumberingLterms} and $\ev_k$ be as in \th\ref{H10Pev}. Let $f$ be as in the definition of the positive answer to Hilbert's tenth problem. Then by definition of $f$ we have that the characteristic function of $S$ is
		$$\X_S(\vec{a}) = f(l,i_1,\ev_k(\vec{a},J(t_{11})),\ldots,\ev_k(\vec{a},J(t_{ik_i})),\ldots, i_l,\ev_k(\vec{a},J(t_{l1})),\ldots,\ev_k(\vec{a},J(t_{ln_l}))).$$
		Since $f$ and $\ev_k$ are recursive functions, we get that $\X_S$ is a recursive function, so $S$ is a recursive set.
	\end{proof} 
	
	\begin{cor} \th\label{H10PSdiophnotrec=>undecidable} 
		Let $(R,\L)$ be a set with a countable language such that $S_c'(\L) = R$. If under all recursive representations $j: R \ra \N$ we have that there exists a Diophantine set $S_j$ such that $S_j$ is not a recursive set under $j$, then Hilbert's tenth problem over $R$ has a negative answer. \qed 
	\end{cor} 
	
	\begin{lemma} \th\label{H10PSdioph=>Srecalgemeen} 
		Let $(R,\L)$ be a set with a countable language. Let $R_0$ be a recursive ring with representation $j$. Suppose that $S_c(R) \subset R_0 \subset S_c'(R)$. If Hilbert's tenth problem over $R$ has a positive answer (under the presentation $j$), then all sets $S \subset R_0^k$ that are Diophantine as a subset of $R^k$ are recursive sets. \qed 
	\end{lemma}

	\subsection{Transferring a Negative Answer to Hilbert's Tenth Problem} 
	
	\begin{lemma}[\th\ref{introH10Pnegtransf}] \th\label{H10Pddiophrec+R1und=>R2und}
		Let $(R_1,\L_1)$ and $(R_2,\L_2)$ be sets with a recursive interpretation of a language. Let $d:R_1 \ra R_2$ be an injective Diophantine and recursive map. If $S \subset R_1$ is a Diophantine set that is not a recursive set, then $d(S)$ is also a Diophantine set that is not a recursive set. If furthermore $d(S) \subset S_c'(\L_2)$, then Hilbert's tenth problem over $(R_2,\L_2)$ has a negative answer. \qed 
	\end{lemma} 
	
	\begin{lemma} \th\label{H10Pddiophequiv+rec+R1und=>R2und} 
		Let $(R_1,\L_1)$ and $(R_2,\L_2)$ be sets with a recursive interpretation of a language. Let $\sim$ be an equivalence relation on $R_2$ and assume that $R_2\modsim$ has a recursive representation such that the projection $\pi: R_2 \ra R_2\modsim$ is a recursive function. Let $d: R_1 \ra R_2\modsim$ be an injective recursive function and a Diophantine equivalence map. If $S\subset R_1$ is a Diophantine set that is not recursive, then $\ol{d}(S) \subset R_2$ is a Diophantine set that is not recursive and Hilbert's tenth problem over $(R_2,\L_2)$ has a negative answer. \qed 
	\end{lemma}  
	
	\subsection{Effective Diophantine maps}\phantom{=}\\
	In order to also transfer positive answers to Hilbert's tenth problem, we need a way to convert Diophantine statements from one ring to another. We formalize this in the following way:
	
	\begin{defi} \th\label{H10Pdefeffdioph}
		Let $(R_1,\L_1)$ and $(R_2,\L_2)$ be sets with a countable language. We take $J_i$ to be a numbering of the $\L_i$-terms of \th\ref{H10PnumberingLterms}.\\
		We say that a Diophantine map $d:(R_1,\L_1)\ra (R_2,\L_2)$ is \textbf{effective} if there exists a corresponding recursive function $D: \bigcup_{n=0}^\infty \N^n \ra \bigcup_{n=0}^\infty \N^n$ such that 
		$$D(l,i_1,J_1(t_{11}),\ldots,J_1(t_{1i_1}),\ldots ,i_l,J_1(t_{l1}),\ldots,J_1(t_{li_l}))$$ 
		corresponds with a system
		$$(s_{11},\ldots,s_{1j_1}) \in S_{j1}',\ldots,(s_{k1},\ldots,s_{kj_m}) \in S_{jm}',$$
		of Diophantine equations in $\L_2$, which means that 
		\begin{align*} &\phantom{=}(k,j_1,J_2(s_11),\ldots,J_2(s_{1j_1}),\ldots,j_k,J_2(s_{k1}),\ldots,J_2(s_{kj_m}))\\
			&= D(l,i_1,J_1(t_{11}),\ldots,J_1(t_{1i_1}),\ldots ,i_l,J_1(t_{l1}),\ldots,J_1(t_{li_l})).
		\end{align*} 
		This system is such that
		\begin{enumerate} 
			\item if $\vec{x}$ is a solution to 
			$$(t_{11},\ldots,t_{1i_1}) \in S_{i_1},\ldots,(t_{l1},\ldots,t_{li_l}) \in S_{i_l}$$
			then there exists vector $\vec{y}$ such that $(d(\vec{x}),
			\vec{y})$ is a solution to 
			$$(s_{11},\ldots,s_{ij_1}) \in S_{j1}',\ldots,(s_{k1},\ldots,s_{kj_m}) \in S_{jm}'.$$
			\item if $(\vec{z},\vec{y})$ is a solution to 
			$$(s_{11},\ldots,s_{ij_1}) \in S_{j1}',\ldots,(s_{k1},\ldots,s_{kj_m}) \in S_{jm}'$$
			then there exists a vector $\vec{x}$ such that $\vec{z} = d(\vec{x})$ and $\vec{x}$ is a solution to 
			$$(t_{11},\ldots,t_{1i_1}) \in S_{i_1},\ldots,(t_{l1},\ldots,t_{li_l}) \in S_{i_l}.$$
		\end{enumerate} 
	\end{defi} 
	
	\begin{thm} \th\label{H10PR2deci+deff=>R1deci}
		Let $(R_1,\L_1)$ and $(R_2,\L_2)$ be sets with a countable language. Let $R_i'$ be a recursive ring such that $S_c(\L_i) \subset R_i' \subset S_c'(\L_i)$ and let $J_i$ be the numbering of the $\L_i$-terms from \th\ref{H10PnumberingLterms}. Suppose that Hilbert's tenth problem over $(R_2,\L_2)$ has a positive answer and $d:R_1 \ra R_2$ is an effective Diophantine map. Then Hilbert's tenth problem over $(R_1,\L_1)$ has a positive answer.
	\end{thm} 
	
	\begin{proof}
		Let $D$ be the function on systems of equations corresponding to $d$. Let $f: \bigcup_{n=0}^\infty \N^n \ra \{0,1\}$ be the function that decides Hilbert's tenth problem over $(R_2,\L_2)$. Then $f\circ D: \bigcup_{n=0}^\infty  \N^n \ra \{0,1\}$ decides Hilbert's tenth problem over $(R_1,\L_1)$, since a system $S$ has a solution $\vec{x}$ if and only if $D(S)$ has a solution, namely $(d(\vec{x}),\vec{y})$. 
	\end{proof} 
	
	\begin{cor}
		Suppose that $(R_1,\L_1)$ and $(R_2,\L_2)$ are sets with a countable language and let $R_1'$ and $R_2'$ be recursive rings such that $S_c(\L_i) \subset R_i' \subset S_c'(\L_i)$. Take $J_i$ to be the numbering of the $\L_i$ terms from \th\ref{H10PnumberingLterms}. If Hilbert's tenth problem over $(R_1,\L_1)$ has a negative answer and $d: R_1 \ra R_2$ is an effective Diophantine map, then Hilbert's tenth problem over $(R_2,\L_2)$ has a negative answer. 
	\end{cor} 
	
	\begin{assumption}
		Let $(R,\L)$ be a set with a language. Using elements of $S_c'(\L)$ as constants does not change which sets are Diophantine, since by definition we have that $\{x\}$ is a Diophantine set for each $x \in S_c'(\L)$. However, we lose information by doing this since we do not know the Diophantine formulas used to describe $\{x\}$. To solve this problem, we assume that only the constants of $\L$, so elements of the set $S_c(\L)$, are used in Diophantine formulas.
	\end{assumption}
	
	\begin{thm}[\th\ref{introEffDiophMap}] \th\label{H10PeffectiveDioph}
		Take $(R_1,\L_1)$ and $(R_2,\L_2)$ sets with a finite language and take $R_i'$ recursive rings such that $S_c(\L_i) \subset R_i' \subset S_c'(\L_i)$. Let $d: (R_1,\L_1) \ra (R_2,\L_2)$ be an injective map. 
		Suppose that for all of the following sets:
		\begin{itemize}
			\item For all constants $c \in \L_1$, the set $\{d(c)\}$;
			\item For all functions $f: R_1^k \ra R_1$ of $\L_1$ the set $\{(d(x_1),\ldots,d(x_k),d(f(x_1,\ldots,x_k))) \mid \vec{x} \in R_1^k\}$;
			\item For all relations $S \subset R_1^k$ of $\L_1$ the set $d(S)$;
			\item The set $d(R_1)$.
		\end{itemize}
		an explicit description of the form 
		$$\{\vec{x} \in R_2^k \mid \ex \vec{y} \in R_2^l,\ \va 1\leq i \leq n,\ (t_{i1}(\vec{x},\vec{y}),\ldots,t_{in_i}(\vec{x},\vec{y})) \in S_i\}$$
		with $t_{ij}$ terms of $\L_2$ and $S_i$ relations of $\L_2$ is known. Then the map $d$ is an effective Diophantine map.
	\end{thm} 
	
	The reason that this theorem holds is that the proof of \th\ref{diophinjdefLwithR} is constructive. To formalize this we need the following lemma.

	\begin{lemma} \th\label{H10PexistencerecD'} 
		\nieuw{Assume that all conditions of \th\ref{H10PeffectiveDioph} hold. Then there exists a recursive function $D':\N \ra \bigcup_{n=0}^\infty \N^n$ such that for an $\L$-term $t$ with $k$ variables the system corresponding to $D'(J_1(t))$ has a set of solutions $S$ such that the projection of $S$ onto the first $k+1$ coordinates gives the set $\{(d(\vec{x}),d(t(\vec{x})))\mid \vec{x} \in R_1^k\}$.} 
	\end{lemma} 
	
	\begin{proof}
		We will define $D'$ by induction on the number of symbols. We will use the function $u$ of \th\ref{H10Puitpakkenterms} for the ring $R_1$.  Since $\L_1$ is finite, the function $F$ of \th\ref{H10Puitpakkenterms} has a finite domain, so it is recursive.\\
		Let $n$ be the input. We do a case distinction.
		\begin{itemize}
			\item If $u(n)_1 = 2^k$ and $F(\rho(1,u(n)_1)) = 0$, we have that $J_1^{-1}(n)$ starts with a constant, so it is a constant. We output the number of the system of equations corresponding to $x_1=d(\rho(1,u(n)_1))$. Since there are only finitely many constants and we know the systems corresponding to $d(c)$ for all $c \in \L$, this is recursive.
			\item If $u(n)_1 = 3^k$, we have that $J_1^{-1}(n)$ starts with a variable $x$, so it a variable $x$. We output the number corresponding to the system $x_1=x_2,\ x_1 \in d(R_1)$.  
			\item If $u(n)_1 = 2^k$ and $F(\rho(1,u(n)_1)) = r>0$, we have that $J_1^{-1}(n)$ starts with an $r$-ary function $f$ applied to $(u(n)_2,\ldots,u(n)_{r+1})$. By induction we know $D'(u(n)_i)$ for $2\leq i \leq r+1$. This means we know the systems $S_2(\vec{x},y,\vec{a}),\ldots,S_{r+1}(\vec{x},y,\vec{a})$ such that the set of solutions of $S_i$ projects onto the set $\{(d(\vec{x}),d(u(n)_i(\vec{x})))\mid \vec{x} \in R_1^{k_i}\}$. By the proof of \th\ref{diophinjdef} we have that
			\begin{align*} 
				\{(d(\vec{x}),d(t(\vec{x})))\mid \vec{x}\in R_1^k\} = \{(d(\vec{x}),d(f(u(n)_2(\vec{x}),\ldots,u(n)_{r+1}(\vec{x}))))\mid \vec{x} \in R_1^k\} \\
				= \{(\vec{u},v)\in R_2^{k+1}\mid \ex (\vec{u_i},v_i)\in \{(d(\vec{x}),d(u(n)_{i+1}(\vec{x})))\mid \vec{x} \in R_1^k\},\ \va 2\leq i \leq r+1,\ \vec{u} = \vec{u_i},\\ 
				(v_2,\ldots,v_{r+1},v) \in \{(d(\vec{y}),d(f(\vec{y})))\mid \vec{y} \in R_1^{r}\}\}.
			\end{align*}
			Let $S(\vec{x},z,\vec{y})$ be the system corresponding to $\{(d(\vec{x}),d(f(\vec{x})))\mid \vec{x} \in R_1^r\}$. Then the system corresponding to $D'(n)$ should be
			\begin{align*}
				\begin{cases} &\va\ 2 \leq i \leq r+1,\ S_i(\vec{u_i},v_i,\vec{a_i})\\
					&\va\ 2\leq i \leq r+1,\ \vec{u} = \vec{u_i}\\
					&S(\vec{v},v,\vec{b}). \end{cases}
			\end{align*}  
			Since $\L$ has only finitely many functions, we can get $S$ recursively. By a variant of \th\ref{H10Pev} we have that substituting variables is recursive. We can get the $S_i$ recursively by inductive definition. Since the codes of systems of equations work in such a way that combining them is recursive, we get that this step in the definition of $D'$ is recursive.
		\end{itemize}
		These are all possibilities since we exclude the use of elements of $R_0$. This means that $D'$ is a recursive function.
	\end{proof} 
	
	\begin{proof}[Proof of \th\ref{H10PeffectiveDioph}]
		By \th\ref{diophinjdef} $d$ is a Diophantine map. We will construct the function $D$ using the proof of \th\ref{diophinjdef} and the previous lemma's. \\
		Suppose that the code of the system
		$$\va 1\leq i \leq n,\ (t_{i1}(\vec{x}),\ldots,t_{im_i}(\vec{x})) \in S_i,$$
		with $t_{ij}$ terms of $\L_1$ and $S_i$ relations of $\L_1$ is the input.\\
		We know from the proof of \th\ref{diophinjdef} that $\vec{x}$ is a solution to this system if and only if there exists a vector $\vec{z}$ such that $(d(\vec{x}),\vec{z}) = (\vec{u},\vec{z})$ is a solution to the system
		\begin{align*}
			\begin{cases} &\va 1\leq i \leq n,\ \va 1\leq j \leq m_i,\ z_{ij} \in d(R_1) \\
				&\va 1 \leq i \leq n,\ \va 1 \leq j \leq m_i,\ (\vec{u},z_{ij}) \in \{(d(\vec{x}),d(t_{ij}(\vec{x})))\} \\
				&\va 1\leq i \leq n,\ \vec{z_i} \in d(S_i). \end{cases} 
		\end{align*}
		We have the recursive function $D'$ from \th\ref{H10PexistencerecD'}, which gives for each $t_{ij}(\vec{x})$ a system $T_{ij}(\vec{x},u,\vec{y_{ij}})$ of Diophantine equations in $\L_2$ such that $T_{ij}$ has the set $\{(d(\vec{x}),d(t_{ij}(\vec{x})),\vec{y_{ij}})\mid \vec{x} \in R_1^k\} $ as set of solutions. Let $S_i'(\vec{z})$ be the system of Diophantine equations corresponding to $\vec{z} \in d(S_i)$ and let $R(z)$ be the system corresponding to $z \in d(R_1)$. We know these systems by the conditions of the theorem. \\
		We define $D$ to map the system $\va 1\leq i \leq n,\ (t_{i1}(\vec{x}),\ldots,t_{im_i}(\vec{x})) \in S_i$ that we started with, to the system $S(\vec{u},\vec{y},\vec{z})$ given by
		\begin{align*}
			\begin{cases}&\va 1\leq i \leq n,\ \va 1\leq j \leq m_i,\ R(z_{ij})\\
				&\va 1\leq i \leq n,\ \va 1\leq j \leq m_i,\ T_{ij}(\vec{u},z_{ij},\vec{y_{ij}}) \\
				&\va 1\leq i \leq n,\ S_i'(\vec{z_i}). \end{cases} 
		\end{align*} 
		We are given $R,\ T_{ij}$ and $S_i'$ recursively. Substitution of variables and combining systems is recursive. Together this implies that $D$ is a recursive function. We have by construction that $\vec{x}$ is a solution to $\va 1\leq i \leq n,\ (t_{i1}(\vec{x}),\ldots,t_{im_i}(\vec{x})) \in S_i$ if and only if $(d(\vec{x}),\vec{y},\vec{z}) = (\vec{u},\vec{y},\vec{z})$ is a solution to $S(\vec{u},\vec{y},\vec{z})$ \nieuw{and all solutions of $S(\vec{u},\vec{y},\vec{z})$ are of the form $\vec{u} = d(\vec{x})$}. This is exactly the requirement for $D$, so $d$ is an effective Diophantine map.
	\end{proof} 
	
	\subsection{Effective Diophantine Equivalence Maps}
	
	\begin{defi} \th\label{H10Pdefeffdiophequiv}
		Let $(R_1,\L_1)$ and $(R_2,\L_2)$ be sets with a countable language. Let $\sim$ be an equivalence relation on $R_2$. We take $J_i$ to be a numbering of the $\L_i$-terms of \th\ref{H10PnumberingLterms}.\\
		We say that a Diophantine equivalence map $d:R_1\ra R_2\modsim$ is \textbf{effective} if there exists a recursive function $D: \bigcup_{n=0}^\infty \N^n \ra \bigcup_{n=0}^\infty \N^n$ such that 
		$$D(l,i_1,J_1(t_{11}),\ldots,J(t_{1i_1}),\ldots ,i_l,J(t_{l1}),\ldots,J_1(t_{li_l}))$$ 
		corresponds to a system 
		$$(s_{11},\ldots,s_{ij_1}) \in S_{j1}',\ldots,(s_{k1},\ldots,s_{kj_m}) \in S_{jm}'$$
		of Diophantine equations in $\L_2$ such that
		\begin{enumerate}
			\item if $\vec{x}$ is a solution to 
			$$(t_{11},\ldots,t_{1i_1}) \in S_{i_1},\ldots,(t_{l1},\ldots,t_{li_l}) \in S_{i_l},$$
			then there exists a vector $\vec{z} \in \ol{d}(\vec{x})$ and a vector $\vec{y}$ such that $(\vec{z},
			\vec{y})$ is a solution to 
			$$(s_{11},\ldots,s_{ij_1}) \in S_{j1}',\ldots,(s_{k1},\ldots,s_{kj_m}) \in S_{jm}'.$$
			\item if $(\vec{z},\vec{y})$ is a solution to 
			$$(s_{11},\ldots,s_{ij_1}) \in S_{j1}',\ldots,(s_{k1},\ldots,s_{kj_m}) \in S_{jm}',$$
			then there exists a vector $\vec{x}$ such that $\vec{z} \in \ol{d}(\vec{x})$ and $\vec{x}$ is a solution to 
			$$(t_{11},\ldots,t_{1i_1}) \in S_{i_1},\ldots,(t_{l1},\ldots,t_{li_l}) \in S_{i_l}.$$
		\end{enumerate} 
	\end{defi} 
	
	\begin{thm} \th\label{H10PR2deci+deffequiv=>R1deci}
		Let $(R_1,\L_1)$ and $(R_2,\L_2)$ be sets with a countable language and let $\sim$ be an equivalence relation on $R_2$. Let $R_i'$ be a recursive ring such that $S_c(\L_i) \subset R_i' \subset S_c'(\L_i)$ and let $J_i$ be the numbering of the $\L_i$-terms from \th\ref{H10PnumberingLterms}. Suppose that Hilbert's tenth problem over $(R_2,\L_2)$ has a positive answer and $d:R_1 \ra R_2\modsim$ is an effective Diophantine equivalence map. Then Hilbert's tenth problem over $(R_1,\L_1)$ has a positive answer.
	\end{thm} 
	
	\begin{proof}
		Let $D$ be the function on systems of equations corresponding to $d$. Let $f: \bigcup_{n=0}^\infty \N^n \ra \{0,1\}$ be the function that decides Hilbert's tenth problem over $(R_2,\L_2)$. \\
		Let $T$ be a system of equations for $R_1$ and let $D(T)$ be the corresponding system in $R_2$. If $T$ has a solution $\vec{x}$, then $D(T)$ has at least one solution, since there exists a $\vec{z} \in \ol{d}(\vec{x})$ and an $\vec{y}$ such that $(\vec{z},\vec{y})$ is a solution to $D(T)$. Suppose now that $D(T)$ has a solution $(\vec{z},\vec{y})$. By the requirements on $D$ there exists a vector $\vec{x}$ such that $\vec{z} \in \ol{d}(\vec{x})$ and $\vec{x}$ is a solution to $T$. This gives that $T$ has a solution if and only if $D(T)$ has a solution. This implies that $f\circ D: \bigcup_{n=0}^\infty  \N^n \ra \{0,1\}$ decides Hilbert's tenth problem over $(R_1,\L_1)$. 
	\end{proof} 
	
	\begin{cor}
		Let $(R_1,\L_1)$ and $(R_2,\L_2)$ be sets with a countable language and let $\sim$ be an equivalence relation on $R_2$. Take $R_1'$ and $R_2'$ be recursive rings such that $S_c(\L_i) \subset R_i'\subset S_c'(\L_i)$. If Hilbert's tenth problem over $(R_1,\L_1)$ has a negative answer and $d: R_1 \ra R_2\modsim$ is an effective Diophantine equivalence map, then Hilbert's tenth problem over $(R_2,\L_2)$ has a negative answer. \qed 
	\end{cor} 
	
	\begin{thm} \th\label{H10PeffectiveDiophEquiv}
		Take $(R_1,\L_1)$ and $(R_2,\L_2)$ sets with a finite language and take $R_i'$ recursive rings such that $S_c(\L_i) \subset R_i' \subset S_c'(\L_i)$. Let $\sim$ be an equivalence relation on $R_2$ and let $d: R_1 \ra R_2\modsim$ be an injective map. 
		Suppose that for the following sets:
		\begin{itemize}
			\item For all constants $c \in \L_1$, the set $\ol{d}(c)$;
			\item For all functions $f: R_1^k \ra R_1$ of $\L_1$ the set $\ol{d}(\{(\vec{x},f(\vec{x}))\mid \vec{x} \in R_1^k\})$;
			\item For all relations $S \subset R_1^k$ of $\L_1$ the set $\ol{d}(S)$;
			\item The set $\ol{d}(R_1)$;
			\item The set $\{(x,y) \in R_2^2\mid x \sim  y\}$ or the set $\{(x,y)\in R_2^2 \mid x,y \in \ol{d}(R_1),\ x\sim y\}$,
		\end{itemize}
		an explicit description of the form 
		$$\{\vec{x} \in R_2^k \mid \ex \vec{y} \in R_2^l,\ \va 1\leq i \leq n,\ (t_{i1}(\vec{x},\vec{y}),\ldots,t_{in_i}(\vec{x},\vec{y})) \in S_i\}$$
		with $t_{ij}$ terms of $\L_2$ and $S_i$ relations of $\L_2$ is known. Then the map $d$ is an effective Diophantine equivalence map.
	\end{thm} 
	
	For the proof we need a variant of \th\ref{H10PexistencerecD'}. 
	
	\begin{lemma} \th\label{H10PexistencerecD'equiv} 
		Suppose that all conditions of \th\ref{H10PeffectiveDiophEquiv} hold. Then there exists a recursive function $D':\N \ra \bigcup_{n=0}^\infty \N^n$ such that for an $\L$-term $t$ with $k$ variables the system corresponding to $D'(J_1(t))$ has a set of solutions $S$ such that the projection of $S$ onto the first $k+1$ coordinates gives the set $\{(\ol{d}(\vec{x}),\ol{d}(t(\vec{x})))\mid \vec{x} \in R_1^k\}$.
	\end{lemma} 
	
	\begin{proof}
		We can nearly copy the proof of \th\ref{H10PexistencerecD'}. We change all occurrences of $d(S)$ for some set $S$ to $\ol{d}(S)$. In the induction basis case that $t$ is a variable we change that we do not use the system $x_1=x_2,\ x_1 \in d(R_1)$, but $x_1 \sim x_2,\ x_1 \in \ol{d}(R_1)$. By the extra assumption in \th\ref{H10PeffectiveDiophEquiv}, we can get a Diophantine form of this system recursively. Instead of \th\ref{diophinjdef} now \th\ref{H10Pdinjdiophequiv} gives us that the proof works. 
	\end{proof} 
	
	\begin{proof}[Proof of \th\ref{H10PR2deci+deffequiv=>R1deci}]
		We have that the second part of the proof of \th\ref{H10Pdinjdiophequiv} goes the same as the second part of the proof of \th\ref{diophinjdef}. The differences in the first part are covered by the use of \th\ref{H10PexistencerecD'equiv} instead of \th\ref{H10PexistencerecD'}. This means we can copy the proof of \th\ref{H10PR2deci+deff=>R1deci}.
	\end{proof} 
	
	\subsection{Transferring MRDP}\phantom{=}\\ 
	In this paragraph, we collect some results on how to use Diophantine and recursive maps to transfer that MRDP holds between rings.
	
	\begin{thm} \th\label{H10PtranferMRDP} 
		Let $(R_1,\L_1)$ and $(R_2,\L_2)$ be recursive rings with a recursive interpretation of a language. Suppose that $d:R_1 \ra R_2$ is a bijective Diophantine and recursive map. If MRDP holds in $R_1$, then MRDP holds in $R_2$. \qed 
	\end{thm} 
	
	\begin{thm} \th\label{H10PMRDPR=>MRDPR[a]} 
		Let $R$ be a recursive commutative domain with a language $\L$. Take $\a$ an element of the algebraic closure of $\Frac(R)$ that is integral over $R$. Let $\L_\a = \L \cup \{\a\}$ be the language on $R[\a]$. If $R \subset R[\a]$ is a Diophantine set and MRDP holds in $(R,\L)$, then MRDP holds in $(R[\a],\L_\a)$. 
	\end{thm} 
	
	\begin{proof}
		The inclusion $\i: R \ra R[\a]$ is a Diophantine map by \th\ref{diophinjdef}, since $R$ is a Diophantine set.\\
		Since $\a$ is integral over $R$, there exist a minimal $n\in \Z_{>0}$ such that we can write $\a^n = r(\a)$ for some $r(X) \in R[X]$ of degree at most $n-1$. Using the commutativity of $R[\a]$, this means we can write each element of $R[\a]$ uniquely in the form $r_0+\a r_1+\ldots + \a^{n-1}r_{n-1}$ with the $r_i \in R$. Now define $d: R^n \ra R[\a]$ by $d(r_0,\ldots,r_{n-1}) = r_0+\a r_1+\ldots + \a^{n-1}r_{n-1}$. Then $d$ is a bijection. We have that $d$ is given by a Diophantine expression in $\i$, namely
		$$d(\vec{x}) = y\ \LRa\ \i(x_1)+\a \i(x_2)+\ldots+ \a^{n-1} \i(x_{n}) = y.$$
		With \th\ref{diophd1dioph+d1<Dd2=>d2dioph} and the fact that $\i$ is an injective Diophantine map this gives that $d$ is a Diophantine map. \\
		We can construct a recursive representation of $R[\a]$ using that of $R$ such that $d: R^n \ra R[\a]$ is a recursive map.\\
		\th\ref{recMRDPR=>MRDPRk} gives that MRDP holds in $(R^n,\L^{(n)})$, since MRDP holds in $(R,\L)$. We have proven that $d$ fulfills the conditions of \th\ref{H10PtranferMRDP}, so MRDP holds in $(R[\a],\L_\a)$. 
	\end{proof} 
	
	\bibliographystyle{amsplain} 
	\bibliography{biblio2}
	
\end{document}